\documentclass[12pt,a4paper]{amsart}
\usepackage{amssymb}

\newcommand{\CC}{{\mathbb{C}}}
\newcommand{\FF}{{\mathbb{F}}}
\newcommand{\GG}{{\mathbb{G}}}
\newcommand{\HH}{{\mathbb{H}}}
\newcommand{\NN}{{\mathbb{N}}}
\newcommand{\QQ}{{\mathbb{Q}}}
\newcommand{\TT}{{\mathbb{T}}}
\newcommand{\ZZ}{{\mathbb{Z}}}

\newcommand{\bi} {\mathbf i}
\newcommand{\br} {\mathbf r}

\newcommand{\bu} {\mathbf u}
\newcommand{\bv} {\mathbf v}
\newcommand{\bw} {\mathbf w}
\newcommand{\bx} {\mathbf x}

\newcommand{\bC} {\mathbf C}
\newcommand{\bG} {\mathbf G}
\newcommand{\bL} {\mathbf L}
\newcommand{\bT} {\mathbf T}
\newcommand{\bU} {\mathbf U}
\newcommand{\bW} {\mathbf W}

\newcommand{\cA} {\mathcal A}
\newcommand{\cB} {\mathcal B}
\newcommand{\cC} {\mathcal C}
\newcommand{\cE} {\mathcal E}
\newcommand{\cF} {\mathcal F}
\newcommand{\cH} {\mathcal H}
\newcommand{\cL} {\mathcal L}
\newcommand{\cO} {\mathcal O}
\newcommand{\cY} {\mathcal Y}

\newcommand{\fp} {\mathfrak p}

\newcommand{\fS} {\mathfrak S}

\newcommand{\Hom}{{{\operatorname{Hom}}}}
\newcommand{\End}{{{\operatorname{End}}}}
\newcommand{\Frac}{{{\operatorname{Frac}}}}

\newcommand{\Ind}{{{\operatorname{Ind}}}}
\newcommand{\IBr}{{{\operatorname{IBr}}}}
\newcommand{\Irr}{{{\operatorname{Irr}}}}

\newcommand{\Mat}{{{\operatorname{Mat}}}}
\newcommand{\op}{{{\operatorname{op}}}}

\newcommand{\rk}{{{\operatorname{rk}}}}
\newcommand{\spets}{{{\operatorname{s}}}}
\newcommand{\sq}{{{\spets,q}}}

\newcommand{\Uch}{{{\operatorname{Uch}}}}
\newcommand{\GL}{\operatorname{GL}}

\newcommand{\tw}[1]{{}^{#1}\!}
\newcommand{\Yop}{{Y^\op}}

\newcommand{\tu}{{\tilde u}}
\newcommand{\hA}{{\hat A}}
\newcommand{\hB}{{\widehat B}}
\newcommand{\tB}{{\widetilde B}}
\newcommand{\tK} {{\tilde K}}

\let\al=\alpha
\let\ze=\zeta
\let\vhi=\varphi
\let\th=\theta
\let\la=\lambda
\let\harr=\hookrightarrow

\newtheorem{thm}{Theorem}[section]
\newtheorem{lem}[thm]{Lemma}

\newtheorem{prop}[thm]{Proposition}
\newtheorem{conj}[thm]{Conjecture}

\newtheorem{ques}[thm]{Question}

\newtheorem{thmA}{Theorem}
\newtheorem{corA}[thmA]{Corollary}
\newtheorem{conjA}{Conjecture}

\theoremstyle{definition}
\newtheorem{rem}[thm]{Remark}
\newtheorem{defn}[thm]{Definition}
\newtheorem{exmp}[thm]{Example}

\begin{document}

\title[The principal block of a $\ZZ_\ell$-spets]{The principal block of a $\ZZ_\ell$-spets\\ and Yokonuma type algebras}

\author{Radha Kessar}
\address{Department of Mathematics, City, University of London EC1V 0HB,
  United Kingdom}
\email{radha.kessar.1@city.ac.uk}

\author{Gunter Malle}
\address{FB Mathematik, TU Kaiserslautern, Postfach 3049,
  67653 Kaisers\-lautern, Germany.}
\email{malle@mathematik.uni-kl.de}

\author{Jason Semeraro}
\address{Heilbronn Institute for Mathematical Research, Department of
  Mathematics, University of Leicester, United Kingdom}
\email{jpgs1@leicester.ac.uk}

\thanks{The first author gratefully acknowledges support from EPSRC grant
 EP/T004592/1. The second author gratefully acknowledges support by the DFG
 --- Project-ID 286237555 -- TRR 195.}

\begin{abstract}
We formulate conjectures concerning the dimension of the principal block of a
$\ZZ_\ell$-spets (as defined in our earlier paper), motivated by analogous
statements for finite groups. We show that these conjectures hold in certain
situations. For this we introduce and study a Yokonuma type algebra for
torus normalisers in $\ell$-compact groups which may be of independent interest.
\end{abstract}

\keywords{$\ell$-adic reflection groups, Yokonuma type algebra, principal
  block, spetses}

\subjclass[2010]{20C08, 20C20, 20F55, 20G40, 16G30; secondary: 20D20, 55R35}

\date{\today}

\maketitle


\section{Introduction}
Let $\ell > 0$ be prime. Motivated by a desire to generalise various
conjectures from finite group theory, in \cite{KMS}, we introduced the notion
of $\ZZ_\ell$-spets, an $\ell$-adic version of the spetses first proposed by
Brou\'e--Malle--Michel \cite{BMM99}. Formally, a $\ZZ_\ell$-spets
$\GG=(W\vhi^{-1},L)$ is a spetsial $\ell$-adic reflection group $W$ on a
$\ZZ_\ell$-lattice $L$ together with an element $\vhi\in N_{\GL(L)}(W)$. Via
the theory of $\ell$-compact groups if $q$ is a power of a prime different
from~$\ell$, then under certain conditions, by \cite{BM07} one can associate a
fusion system $\cF$ on a finite $\ell$-group~$S$ to the pair $\GG(q):=(\GG,q)$.
In this situation, $S$ should be considered as a Sylow $\ell$-subgroup of
$\GG(q)$. Using $\cF$, we also attached in \cite{KMS} a ``principal
$\ell$-block'' $B_0$ to $\GG(q)$, with defect group $S$ and associated set
$\Irr(B_0)$ of characters defined in terms of a collection of Hecke algebras
associated to the centralisers of elements of $S$ in $\cF$ as explained in
\cite{KMS}. The construction
of $B_0$ is modelled on the case of blocks of finite groups of Lie type in
non-describing characteristic. When $W$ is rational and under some natural
assumptions on $\ell$, we recover the $\ell$-fusion system $\cF_\ell(\GG(q))$
and principal $\ell$-block $B_0(\GG(q))$ of the associated finite reductive
group $\GG(q)$. See Section~\ref{sec:sec3} for a description of $B_0$ and
comparison with the rational case.

A primary concern in \cite{KMS} was the translation of certain local-global
statements in modular representation theory to purely local statements using
the language of $\ZZ_\ell$-spetses. Here, our considerations are more on the
global side, so that actual degrees of characters in $\Irr(B_0)$ (as defined
in \cite[Def.~6.7]{KMS}) play a more significant role. Firstly, we investigate
the dimension
$$\dim(B_0):=\sum_{\gamma \in \Irr(B_0)} \gamma(1)^2\in\ZZ[x]$$
of $B_0$. Motivated by results concerning the divisibility properties of this
number for principal blocks of finite groups, we make the following conjecture:

\begin{conjA}   \label{c:main}
 Let $\GG=(W\vhi^{-1},L)$ be a simply connected $\ZZ_\ell$-spets for which
 $\ell$ is very good, let $q$ be a power of a prime different from $\ell$, and
 $B_0$ be the principal $\ell$-block of $\GG(q)$ with defect group $S$. Then,
 \begin{itemize}
  \item[(1)] $\left(\dim(B_0)|_{x=q}\right)_\ell=|S|$; and
  \item[(2)] $\left(\dim(B_0)|_{x=q}\right)_{\ell'}
    \equiv|W_{\vhi^{-1}\zeta^{-1}}|_{\ell'}\pmod \ell$,
   where $\zeta\in\ZZ_\ell^\times$ is the root of unity with
   $q\equiv\zeta\pmod\ell$, and $W_{\vhi^{-1}\zeta^{-1}}$ is the associated
   relative Weyl group (see \cite[Thms~2.1 and~3.6]{KMS} for the definition and
   role of the relative Weyl group).
 \end{itemize}
\end{conjA}

For a finite group $G$ with principal $\ell$-block $B_0$ the equality
$\left(\dim(B_0)\right)_\ell=|S|$ is closely related to the fact that $|S|$
divides the dimension of each projective indecomposable module $\Phi_\nu$
associated to $\nu\in\IBr(B_0)$. When $B_0$ is as in Conjecture~\ref{c:main},
and under the additional hypotheses that $W$ is an $\ell'$-group, $\varphi $
is trivial and $q\equiv 1\pmod\ell$, in Section~\ref{subsec:deco} we formulate
analogues of $\IBr(B_0)$, decomposition numbers and $\deg \Phi_\nu$ for which
we make the following additional conjecture:

\begin{conjA}   \label{c:main2}
 In the setting of Conjecture~{\rm\ref{c:main}} if $W$ is an $\ell'$-group,
 $\varphi$ is trivial and $q\equiv1\pmod\ell$ we have
 \begin{itemize}
  \item[(3)] $|S|$ divides $\deg\Phi_\nu|_{x=q}$ for all $\nu\in\IBr(B_0)$.
 \end{itemize}
\end{conjA}

Conjecture~\ref{c:main} combines the statements of Conjectures~\ref{conj:dim B}
and \ref{conj:quot} and Conjecture~\ref{c:main2} is restated as
Conjecture~\ref{conj:decomp} in Section~\ref{sec:sec3}. In
Proposition~\ref{prop:rational}, we prove that Conjectures~\ref{c:main}
and~\ref{c:main2} hold when $W$ is rational or when $q\equiv 1\pmod\ell$ and
$W$ is primitive. The general case for
Conjecture~\ref{c:main} appears quite elusive, so we choose to simplify matters
by assuming that $W$ is an $\ell'$-group, $q\equiv 1\pmod\ell$ and $\vhi=1$.
Under these conditions we are able to show that our conjectures are closely
related to certain previously-studied properties of Hecke algebras which we
discuss next.
 
Recall that the generic Hecke algebra $\cH(W,\bu)$ of $W$ is a certain quotient
of the group algebra of the braid group $B(W)$ of $W$ where $\bu$ is a set of
parameters indexed by conjugacy classes of distinguished reflections in $W$. To
each irreducible character $\chi$ of $\cH(W,\bu)$ over a splitting field one
associates a Schur element $f_\chi$ which, in turn, is used to define a
canonical form
$$t_{W,\bu}:=\sum_{\chi \in \Irr(\cH(W,\bu))} \frac{1}{f_\chi} \chi$$
on $\cH(W,\bu)$. This is conjectured to be a symmetrising form \cite{BMM99}
(see Conjecture~\ref{conj: symmform Hecke}). It is expected that the form
$t_{W,\bu}$ behaves well under restriction to parabolic subgroups and with
respect to the natural map $B(W) \rightarrow W$; if that is the case we say
$\cH(W,\bu)$ is \emph{strongly symmetric} (see Definition~\ref{def:strongly
quasi}). Our main result is the following:

\begin{thmA}   \label{t:main}
 Let $\GG=(W,L)$ be a simply connected $\ZZ_\ell$-spets. Suppose that
 $\cH(W,\bu)$ is strongly symmetric, $W$ is an $\ell'$-group and
 $q\equiv1 \pmod\ell$. Then Conjectures~{\rm\ref{c:main}} and~{\rm\ref{c:main2}}
 hold. 
\end{thmA}

Since the strongly symmetric condition is known to hold in many cases (see
Proposition~\ref{p:cases}), we obtain:

\begin{corA}   \label{cor:main}
 Suppose $W$ is a spetsial irreducible $\ell'$-group, $\varphi=1$ and
 $q\equiv1 \pmod\ell$.
 Then Conjecture~{\rm\ref{c:main}} holds for $ \GG= (W,L)$,
 and
 Conjecture~{\rm\ref{c:main2}} also holds, except possibly when $W$ is
 primitive of rank at least~$3$.
\end{corA}

The method behind our proof of Theorem~1 may be even more interesting than the
theorem itself. For our approach we introduce and study a Yokonuma type algebra
$\cY $ attached to an arbitrary finite $\ell$-adic reflection group. We
conjecture that $\cY$ is finitely generated and free over its base ring
(Conjecture~\ref{conj:free}) and show that this holds in many cases. It seems
likely that if $\ell$ is very good for $(W,L)$, then $\cY$ is a symmetric
algebra with respect to a trace form whose Schur elements are derived
from those of various parabolic subalgebras of $\cH(W,\bu)$
(Question~\ref{conj:trace}). 

When $\GG$ is a rational spets which satisfies the hypotheses of
Theorem~\ref{t:main}, the principal block $B_0$ of $\GG(q)$ is known to be
Morita equivalent to $SW$ over a finite extension of $\ZZ_\ell$. The
endomorphism algebra of the permutation module which captures this equivalence
is a Yokonuma algebra isomorphic to a certain $\ZZ_\ell$-algebra specialisation
$\cY_\psi$ of $\cY$ (see Section~\ref{sub:classical yokonuma}). Moreover, the
quantities $\dim(B_0)$ and $\deg \Phi_\nu$ ($\nu\in\IBr(B_0)$) can be
re-expressed in terms of the corresponding Schur elements of $\cY_\psi$.

This motivates the study of $\cY_\psi$ and associated numerical properties for
arbitrary $\ell$-adic reflection groups. In Theorem~\ref{thm:O finite generation} we show that    if    the parabolic subgroups  of the underlying  $\ell$-adic reflection group  are generated by reflections, then $\cY_\psi$ is finitely generated and free over
$\ZZ_\ell$; a result of K\"{u}lshammer--Okuyama--Watanabe then allows us to
conclude that when the order of $W$ is relatively prime to $\ell$, then
$\cY_\psi$ is isomorphic to the group algebra of $SW$ (see
Theorem~\ref{thm:KOW2}). We obtain Theorem~1 by applying
general results concerning the arithmetic behaviour of Schur elements in
symmetric algebras (as discussed in Section~\ref{ss:sym}) to $\cY_\psi$, by
playing off the form on $\cY_\psi$ inherited from $\cY$ against the standard
symmetrising form on the group algebra of $SW$.

We close the introduction by discussing how the hypotheses on $\ell$ and $W$ in
Theorem~\ref{t:main} might be relaxed.
If $\ell \mid |W|$ we
have carried out explicit computations in support of Conjecture~\ref{c:main}
when $q\equiv 1\pmod\ell$, such as when $W=G(e,1,3)$ for $\ell=3$ and when
$W=G_{29}$ for the bad prime $\ell=5$ assuming $\Irr(B_0)$ is defined
appropriately (see \cite[Rem.~6.16]{KMS}). In these situations, $\Irr(B_0)$ is
only partially describable in terms of the Schur elements of the Yokonuma
algebra $\cY_\psi$. Even for Weyl groups it is a major open problem to find
a description involving the Schur elements of a suitable larger algebra.
 
\medskip
\subsection*{Structure of the paper}
In Section~\ref{sec:back} we collect some general material necessary for our
proofs. In particular, we discuss certain divisibility properties of Schur
elements in symmetric algebras. In Section~\ref{sec:Hecke} we recall the
construction of Hecke algebras $\cH(W,\bu)$ for complex reflection groups $W$
and the origin of the trace form $t_{W,\bu}$, define the property of being
strongly symmetric (Definition~\ref{def:strongly quasi}) and introduce the
relevant specialisations. In Section~\ref{sec:sec3}, we recall the description
of the principal block of a $\ZZ_\ell$-spets, introduce the notion of dimension
and formulate our main conjectures. Our new Yokonuma type algebra~$\cY$
attached to an $\ell$-adic reflection group $(W,L)$ is defined and studied in
Section~\ref{sec:Yoku}. In
Theorem~\ref{thm:constr rep} we describe the structure of $\cY$ over the field
of fractions of the base ring and in Theorem~\ref{thm:KOW2} we obtain, under
additional conditions, a similar structural result for $\ell$-adic
specialisations of~$\cY$. 
In Sections~\ref{sub:free} and \ref{subsec:trace form} we formulate general
freeness and symmetrising form conjectures for~$\cY$ and show that the freeness
conjecture holds for Coxeter groups (Theorem~\ref{thm:free Cox}) and for most
of the imprimitive complex reflection groups (Theorem~\ref{thm:free G(e,p,n)}).
We also discuss the relationship of $\cY$ to the algebra considered by Marin
in~\cite{Ma18a,Ma18b} and in Section~\ref{sub:classical yokonuma} we show that
when $(W,L)$ arises from a Weyl group the principal block of the classical
Yokonuma algebra over $\ZZ_\ell$ is a certain specialisation of~$\cY$.
Section~\ref{subsec:pf thm1} contains the proof of Theorem~\ref{t:main}
(Theorem~\ref{thm:omni}) and Corollary~\ref{cor:main}. 

\medskip
\noindent{\bf Acknowledgement:}
We thank Maria Chlouveraki for providing pointers to results in
\cite[\S34]{Lu7}, Markus Linckelmann for useful discussions on several aspects
of the paper and in particular on the proof of Theorem~\ref{thm:KOW2} and
Burkhard K\"ulshammer for providing background and references for
Proposition~\ref{prop:dimgreen}. We are indebted to Ivan Marin for his
pertinent comments on an earlier version. We also thank the referee for
numerous constructive comments and recommendations which have helped to improve
the readability of the paper. 

\section{Background Material}   \label{sec:back}

\subsection{Finite generation of modules}
First we record a few general facts on finite generation of modules. The first
is a variation on Nakayama's lemma for nilpotent ideals which allows for the
dropping of the finite generation hypothesis. Let $R$ be a commutative ring
with $1$. Recall that the Jacobson radical $J(R)$ of $R$ is the intersection
of all maximal left ideals of $R$.

\begin{lem}   \label{lem:Nakayama}
 Let $I \subseteq J(R)$ be an ideal of $R$ and let $M,N$ be $R$-modules with
 $N\subseteq M $ and $M = N+ I M$. Suppose that either $M/N$ is finitely
 generated or that $I$ is nilpotent. Then $N=M$.
\end{lem}

\begin{proof}
The case that $M/N$ is finitely generated is the usual case of the Nakayama
Lemma. Suppose that $I^r =0$ for some $r\ge1$. By hypothesis, $M/N = I(M/N)$.
Hence $M/N = I^r (M/N) =0$, showing that $M=N$.
\end{proof} 

\begin{lem}   \label{lem:lifting}
 Let $ R$ be a discrete valuation ring with
 uniformiser $\pi$. Suppose that $ N \subseteq M $ are $R$-modules such
 that $M/\pi M $ is finitely generated, $ M/N $ is finitely generated and
 torsion free, and $\pi^ r N = 0 $ for some positive integer $r$. Then $M$
 is finitely generated.
\end{lem} 

\begin{proof}
Since $M/N$ is finitely generated, it suffices to show that $N$ is finitely
generated. Now $\pi^r N =0$. So, in order to show that $N$ is finitely
generated it suffices to show that $\pi^iN/\pi^{i+1}N$ is finitely generated
for any $i\geq 1$. Multiplication
by $\pi^i $ induces a surjective $R$-module homomorphism from $N/\pi N$ onto
$\pi^iN/\pi^{i+1}N$, hence it suffices to show that $N/\pi N$ is finitely
generated. By hypothesis, $M/\pi M$ is finitely generated and $N/(\pi M\cap N)$
is isomorphic to a submodule of $M/\pi M$. Thus, as $R$ is Noetherian,
$N/(\pi M \cap N)$ is
finitely generated. But since $M/N$ is torsion free, $\pi M\cap N = \pi N$.
\end{proof} 

\subsection{Clifford theory}
The following is standard Clifford theory adapted to quotients of infinite
group algebras with respect to finite normal subgroups. For a ring $R$ and
ideal $I$, we will regard without further comment an $R/I$-module as an
$R$-module via pullback along the canonical homomorphism $R\to R/I$.

\begin{lem}   \label{lem:infclifford}
 Let $G$ be a group, $T$ a finite normal subgroup of $G$ and $K$ a field of
 characteristic~$0$. For $\th \in\Irr_K(T)$, denote by $e_\th$ the
 corresponding primitive central idempotent of $KT$ and by $G_\th$ the
 stabiliser of $\th$ in $G$ of (finite) index $n_\th:=|G:G_\th|$. Let $I$ be an
 ideal of $KG$ and set $I_\th:= e_\th I e_\th = I \cap e_\th KGe_\th$.
 \begin{enumerate}
  \item[\rm(a)] There is a $K$-algebra isomorphism
   $$ KG /I \cong \prod_\th \Mat_{n_\th} ( e_\th KGe_\th/I_\th) $$ 
   where $\th$ runs over a set of representatives of $G$-orbits on $\Irr_K(T)$.
   Moreover,
   $e_\th KGe_\th/I_\th = e_\th KG_\th e_\th/I_\th = KG_\th e_\th/I_\th$ for all
   $\th \in \Irr_K(T) $.
  \item[\rm(b)] The map $(\th,U)\mapsto \Ind_{G_\th} ^G U$ induces a bijection
   between the set of pairs $(\th,U)$, where $\th$ runs over representatives of
   $G$-orbits on $\Irr_K(T)$ and $U$ runs over a set of isomorphism classes of
   simple $KG_\th e_\th /I_\th $-modules, and the set of isomorphism classes of
   simple $KG/I$-modules.
 \end{enumerate}
\end{lem} 

\begin{proof}
Set $ Y =KG/I $ and for $a \in KG$ denote by $\bar a$ its image in $Y$.
Let $X$ be the set of orbits of the conjugation action of $G$ on $\Irr_K(T)$
and for each $x \in X$, let $e_x = \sum_{\th \in x}e_\th$. 
So $\{e_x\mid x\in X\}$ is a set of pairwise orthogonal central idempotents of
$KG$ with $1_{KG}= \sum_{x\in X} e_x$. Consequently,
$\{\bar e_x\mid x \in X\}$ is a set of pairwise orthogonal central idempotents
of $Y$ with $1_{Y} = \sum_{x\in X}\bar e_x$. Here, we abuse notation to allow
for the possibility that some $\bar e_x$ are equal to zero. Thus
$$ Y= \bigoplus _{x\in X} Y \bar e_x $$
and this is also a decomposition of $Y$ into a direct product of
$K$-algebras. Let $x\in X$ and $\th\in x$. Then,
$$ \bar e_x = \sum _{ g \in G/G_\th } g \bar e_\th g^{-1} $$
is a decomposition of $\bar e_x$ into orthogonal, conjugate idempotents. Hence,
$$ Y \bar e_x \cong \Mat_{n_\th} (\bar e_\th Y \bar e_\th )$$
as $K$-algebras. Now the first assertion of (a) follows since the inclusion of
$e_\th KGe_\th$ in $KG$ induces an isomorphism
$e_\th KG e_\th / I_\th \cong \bar e_\th Y \bar e_\th $. 

Since $e_\th g e_\th =0$ for any $g \in G\setminus G_\th$, and since $e_\th$
is central in $KG_\th $,
$$e_\th KGe_\th /I_\th = e_\th KG_\th e_\th /I_\th = KG_{\th}e_\th/I_\th ,$$
proving the second assertion of (a).

Let $V$ be a simple $Y$-module. There exists a unique $x\in X$ such that
$\bar e_x V\ne 0$ (equivalently $e_x V= V)$. By the statement and proof of~(a),
if $\th \in x$, then 
$$ V = \bigoplus_ {g \in G/G_\th} g \bar e_\th V, \eqno{(*)}$$
$\bar e_\th V$ is a simple $\bar e_\th Y \bar e_\th$-module, and the map
$V\mapsto (\th,\bar e_\th V)$ induces a bijection between the set of
isomorphism classes of simple $Y$-modules and the set of pairs $(\th,U)$,
where $\th$ runs over representatives of $G$-orbits on $\Irr_K(T)$ and $U$ runs
over a set of isomorphism classes of simple
$\bar e_\th Y\bar e_\th$-modules. 
Finally, identifying $\bar e_\th Y \bar e_\th$-modules with
$KG_\th e_\th/I_\th$-modules via the isomorphism
$KG_\th e_\th / I_\th \cong\bar e_\th Y \bar e_\th$ given in the proof
of (a), the equation $(*)$ gives
that $V =\Ind_{G_\th} ^G (e_\th V)$. This proves (b).
\end{proof}

\subsection{Symmetrising forms and divisibility}\label{ss:sym}
Let $R$ be a commutative ring with $1$ and let $Y$ be a symmetric $R$-algebra
which is free and finitely generated as $R$-module. If $X$ is an $R$-basis
of~$Y$, then any symmetrising form $t: Y\to R$ determines a dual basis
$X^\vee=\{x^\vee\ |\ x\in X\}$ satisfying
$$t(xy^\vee)=\begin{cases} 1& \text{for $x=y\in X$}\\
                           0& \text{for $x,y\in X$, $x\neq y$}.\end{cases}$$
The \emph{relative projective element of $Y$ in $Z(Y)$ with respect to $t$}
is defined by $z_t := \sum_{x\in X} xx^\vee$. It depends on~$t$ but not on the
choice of the basis $X$ (see \cite[Sec.~2.11, ~2.16]{Li18} for details). 

Now suppose that $t$ is a symmetrizing form on $Y$ and let $s:Y\to R$ be a
symmetric form (also known as trace form). Then there exists $u\in Z(Y)$ such
that $s=t_u$, where $t_u\in Y^*:=\Hom_R(Y,R)$ is defined by $t_u(x):= t(ux)$
for $x\in Y$. The map $u\mapsto t_u$ is an $R$-module isomorphism between
$Z(Y)$ and the $R$-module of symmetric forms
on $Y$. Moreover, $t_u$ is a symmetrising form if and only if
$u \in Z(Y)^\times$. If $X$ is an $R$-basis of $Y$ and $X^\vee$ is the dual
basis with respect to $t$ and if $u \in Z(Y)^\times$, then the dual basis
of~$X$ with respect to $t_u$ is equal to $X^\vee u^{-1}$, and hence the
relative projective element in $Z(Y)$ with respect to $t_u$ is equal
to~$z_t u^{-1}$.

If $Y$ is a split semisimple algebra over a field $K$ then for any symmetric
form $s: Y \to K$, there exist elements $s_\chi \in K$,
$\chi\in\Irr_K(Y)$, such that $s =\sum_{\chi \in \Irr_K(Y)} s_\chi\chi$.
Moreover, we claim that $s$ is a symmetrising form if and only if all $s_\chi$
are non-zero. For this, note that $s$ is a symmetrising form for $Y$ if and
only of its restriction to any block of $Y$ is a symmetrising form for the
block. Thus we may assume that $Y$ is split simple, that is, a matrix algebra
over $K$ and here the claim is immediate from the fact that the trace map on a
matrix algebra is a symmetrising form (see for instance
\cite[Thm~2.11.3]{Li18}).

\begin{lem} \label{lem:schurproj}
 Suppose that $K$ is a field and $Y$ is a split semisimple $K$-algebra.
 Let $s: Y\to K$ be a symmetrising form and let $s_\chi$, $\chi\in \Irr _K(Y)$,
 be elements of $K$ such that $s =\sum_\chi s_\chi\chi$. Let $z=z_s\in Z(Y)$. 
 \begin{enumerate}
  \item[\rm(a)] If $V$ is a $Y$-module affording the irreducible character
   $\chi$, then the trace of $z^{-1}$ on $V$ equals $s_\chi$.
  \item[\rm(b)] The trace of $z^{-2}$ in the left regular representation of
   $Y$ equals $\sum_{\chi\in\Irr_K(Y)} s_\chi^2$.
 \end{enumerate}
\end{lem}

\begin{proof}
A straightforward calculation using the standard bases of matrix algebras shows
that
$$z = \sum_{\chi\in\Irr_K(Y) } s_\chi^{-1} \chi (1) e_\chi$$
where $e_\chi\in Z(Y)$ is the central idempotent corresponding to $\chi$. 
In particular, $z$ is invertible in $Y$ with
$z^{-1} = \sum_\chi s_\chi \chi (1) ^{-1} e_\chi$.
The trace formula is an immediate consequence of this.

For (b) we note that the $K$-linear map
$$s\otimes s: Y\otimes_K \Yop\to K,\qquad
  (s\otimes s) (y\otimes y'):= s(y)s(y'),$$
is a symmetrising form on $Y\otimes_K\Yop$ with corresponding relative
projective element $z\otimes z$. Further,
$Y \cong \bigoplus_{\chi\in\Irr_K(Y)}V_\chi\otimes V_\chi^*$ as
$(Y\otimes_K \Yop)$-modules, where $V_\chi $ is an irreducible $Y$-module
with character $\chi$ and $V_\chi^*$ is an irreducible $\Yop $-module with
character $\chi$. So, by part~(a) applied with $Y\otimes_K \Yop$ in place of
$Y$, $s\otimes s$ in place of $s$ and $V_\chi\otimes V_\chi^*$ in place of $V$,
$\sum_\chi s_\chi^2$ equals the trace of $z^{-1}\otimes z^{-1}$ on $Y$. Since
$z$ is central in $Y$, this trace is just the trace of left multiplication by
$z^{-2}$ on $Y$.
\end{proof}

\begin{lem}   \label{lem:reversedivisibility}
 Let $G$ be a finite group and $K$ a field such that $KG$ is split semisimple.
 Let $t$ be the canonical symmetrising form on $KG$. Let $u \in Z(KG)^\times$,
 let $\al\in K$ be the coefficient of~$1$ when $u^2$ is written as a
 $K$-linear combination of elements of $G$ and suppose that
 $t_u =\sum_\chi s_\chi\chi$, $s_\chi\in K$. Then
 $$\sum_{\chi\in\Irr_K(Y)} s_\chi ^2 = \frac{\al}{|G|}.$$
\end{lem}

\begin{proof}
By a straightforward calculation using the set of group elements as basis, the
relative projective element of $KG$ with respect to $t$ is $|G| 1_{K G}$.
Hence the relative projective element of $KG$ with respect to $t_u$ is
$|G|u^{-1}$. Now the result follows from Lemma~\ref{lem:schurproj}
since for any $y\in KG$, the trace of the action of $y$ on $KG$ via the left
regular representation equals $|G|\beta$, where $\beta$ is the coefficient of
$1$ in the standard basis presentation of~$y$. 
\end{proof} 

\begin{lem}   \label{lem:otherdivisibility}
 Let $\cO$ be a complete discrete valuation ring with field of fractions
 $K$ of characteristic zero and residue field $\cO/J(\cO)$ of
 characteristic~$\ell$. Let $G =T\rtimes W$ be the semidirect product of a
 finite
 abelian $\ell$-group $T$ with an $\ell'$-group $W$ acting faithfully on $T$.
 Let $t$ be the canonical symmetrising form on $KG$ and $s$ a symmetrising form
 on $KG$ with $s = t_u$, $u\in Z(KG)^\times$. Let $\al$ be the coefficient
 of~$1$ when $u^2$ is expressed as a $K$-linear combination of elements of $G$. 

 Suppose that the restriction of $s$ to $\cO G$ takes values in $\cO$
 and for any $x\in T$, $s(x) =\delta_{x,1}$. Then $u\in Z(\cO G)$, $\al\in\cO$
 and $\al\equiv 1\pmod\ell$.
\end{lem}

\begin{proof}
The restriction $t':= t|_{\cO G}$ is a symmetrising form $\cO G\to\cO$ and by
hypothesis $s':=s|_{\cO G}:\cO G\to\cO$ is a symmetric form. Thus there
exists $u'\in Z(\cO G)$ such that $s' = t'_{u'}$. Extending scalars, by
linearity we have that $s = t_{u'}$. Since by hypothesis $s = t_u$, we conclude
$u =u' \in Z(\cO G)$. This proves the first and second assertions.

For a conjugacy class $C$ of $G$, denote by $\hat C$ the corresponding class
sum. Let
$$u = 1 + \sum_C \al_C \hat C $$
where $C$ runs over the non-identity conjugacy classes of $G$. Then for
$1\ne x\in T$ we have 
$0 = s(x) = t( ux) = \al_D $
where $D$ is the conjugacy class of $G$ containing $x^{-1}$.
It follows that 
$$ u = 1 + \sum_{C\in \cC} \al_C\hat C $$
where $\cC$ is the set of conjugacy classes of $G\setminus T$ and
$$ \al= 1 + \sum_{C\in \cC} \al_C \al_{C^{-1} } | C|, $$
where $C^{-1}$ denotes the class containing the inverses of the elements of~$C$.
Since $W$ acts faithfully on $T$, $\ell$ divides $| C|$ for all
$C\in \cC$. Since $u\in Z(\cO G)$, all $\al_C$ are elements of $\cO$ and
we obtain the last assertion.
\end{proof} 

The following is a consequence of Tate duality for symmetric algebras over
complete discrete valuation rings exhibited in \cite{EGKL}. For an integral
domain $\cO$ with field of fractions $K$ and $Y$ an $\cO$-algebra which is
finitely generated free as $\cO$-module, we set $KY:=K\otimes_\cO Y$. 

\begin{prop}   \label{lem:projschur}
 Suppose $\cO$ is a complete discrete valuation ring with field of fractions
 $K$ of characteristic zero, $Y$ is a symmetric $\cO$-algebra such
 that $ KY$ is split semisimple. Let $s:K Y\to K$ be a symmetrising form with
 $s = \sum_\chi s_\chi\chi$. Suppose that the restriction of $s$ to $Y$ takes
 values in $\cO$.

 Let $U$ be a projective $Y$-module and suppose that there is an isomorphism of
 $KY$-modules
 $$K\otimes_\cO U \cong\bigoplus_{\chi \in \Irr(KY) } V_\chi^{d_\chi} ,$$
 where $V_\chi$ is a simple $KY$-module with character $\chi$ and
 $d_\chi\in\NN_0$. Then
 $$\sum_{\chi\in\Irr(KY)} s_\chi d_\chi \in \cO.$$
\end{prop}

\begin{proof}
Let $t: Y\to\cO$ be a symmetrising form and let $z$ be the relative projective
element of $Y$ with respect to $t$. By \cite[Prop. 2.2]{EGKL}, for any
$\gamma\in\End_Y(U)$, the trace of $z^{-1} \gamma$ on $K\otimes_{\cO} U$ lies
in $\cO$. Here by $z^{-1}\gamma\in\End_{KY}(K \otimes_{\cO} U)$ we denote the
composition of (the extension to $K$ of) $\gamma$ with multiplication by
$z^{-1}$.

Denote by $\tilde t : KY \to K$ the $K$-linear extension of $t$ to $KY$. Then
$\tilde t$ is a symmetrising form for $KY$ with relative projective element
$z$. Hence $s = \tilde t_u$ for some $u \in Z(KY)$. Since the restriction of
$s$ to $Y$ takes values in $\cO$, we have as in the proof of
Lemma~\ref{lem:otherdivisibility} that $u \in Z(Y)$. Let $\gamma: U\to U$ be
multiplication by $u$. Since $u\in Z(Y)$, $\gamma\in \End_Y(U)$.
Thus the trace of $z^{-1} \gamma $ on $K\otimes_{\cO} U$ is an element of~$\cO$.
Further, $z u^{-1}$ is the relative projective element of $KY$ with respect to
$s$ and $z^{-1}\gamma$ is multiplication by $(zu^{-1} )^{-1}$. Thus the result
follows from Lemma~\ref{lem:schurproj}.
\end{proof}

The first statement of the next proposition is a theorem of Brauer. The second
is also well known to experts; we provide a proof of it for the convenience of
the reader.

\begin{prop}   \label{prop:dimgreen}
 Let $G$ be a finite group, $S$ a Sylow $\ell$-subgroup of $G$ and $H$ a
 subgroup of $G$ containing $N_G(S)$. Denote by $B_0(G)$ (respectively
 $B_0(H)$) the principal $\ell$-block of $G$ (respectively $H$). Then,
 $$ (\dim B_0(G))_\ell = (\dim B_0(H) )_\ell = |S| $$
 and 
 $$ (\dim B_0(G))_{\ell'}\equiv(\dim B_0(H))_{\ell'}\pmod\ell.$$
\end{prop} 

\begin{proof}
Let $\cO$ be a complete discrete valuation ring in characteristic zero with
algebraically closed residue field of characteristic $\ell$. Let $b\in\cO G$
be the block idempotent corresponding to $B_0(G)$ and $c\in\cO H$ the one
corresponding to $B_0(H)$. The group $S$ is a defect group of the principal
block $B_0(G)$. Hence as $\cO[G\times G]$-module, $B_0(G)$ has vertex
$\Delta S=\{(x,x)\mid x\in S\}$ (see \cite[Rem.~6.7.14]{Li18}). Further, the
$\cO[H \times H]$-module $B_0(H)$ is the Green correspondent of $B_0(G)$ in
$H \times H$ (see \cite[Thm~6.7.2]{Li18}, and note that by
\cite[Thm~6.13.14]{Li18}, $B_0(G)$ and $B_0(H)$ are Brauer correspondents).
Thus, by properties of the Green correspondence (see \cite[Thm~5.2.1]{Li18})
$$\Ind_{H \times H} ^{G \times G} (\cO H c ) = \cO Gb \oplus Y$$
where every indecomposable $\cO[G \times G]$-module summand of $Y$ has vertex
of order strictly smaller than $|\Delta S| = |S|$. Comparing $\cO$-ranks, 
$$\frac{|G|^2} {|H|^2} \rk (\cO Hc)= \rk (\cO Gb) + \rk (Y).$$
By a standard application of Green's indecomposability theorem and the
properties of vertices and sources (see \cite[Sec.~5.1, Thm~5.12.3]{Li18}),
the $\ell$-part of the rank of any indecomposable $\cO X$-lattice, for $X$ a
finite group, is greater than or equal to $\frac{|X|_\ell}{|Q|}$ where $Q$ is
a vertex of the lattice. Thus,
$(\rk(Y))_\ell$ is strictly greater than $|S|$ and we obtain
$$\frac{|G|^2}{|H|^2} \rk(\cO Hc)\equiv \rk(\cO Gb)\pmod {\ell|S|}.$$
Now by a theorem of Brauer (see \cite[Thm~6.7.13]{Li18}), the $\ell$-part of
$\dim( B_0(G)) =\rk(\cO  Gb) $ equals $|S|$ and similarly for $B_0(H)$. The result follows
since $|G:H|\equiv 1\pmod\ell$ by Sylow's theorem.
\end{proof}

\subsection{Tits deformation theorem}
We recall some features of Tits' deformation theorem. Let $R$ and $R'$ be
integral domains with field of fractions $K$ and $K'$ respectively and let
$\psi:R\to R'$ be a ring homomorphism. Let $Y$ be an $R$-algebra which is
finitely generated and free as $R$-module and let $Y'$ denote the
$R'$-algebra $R'\otimes Y$ obtained via extension of scalars through $\psi$.
Let $t:Y\to R$ be an $R$-linear map and let $t':Y'\to R'$ be its
$R'$-linear extension through $\psi$.

\begin{thm}   \label{t:Tits}
 Suppose that $KY$ and $K'Y'$ are both split semisimple.
 \begin{enumerate}
  \item[\rm(a)] The map $\psi$ induces a bijection $\Irr(KY)\to\Irr(K'Y')$,
   $\chi\mapsto\chi'$, such that
   $$\chi'(1\otimes y) = \psi^*(\chi(y))\qquad\text{for all $y \in Y$} .$$ 
   Here $\psi^*: R^*\to K^*$ is an extension of $\psi$ to the integral closure
   $R^*$ of $R$ in $K $ and $K^*$ is some extension field of $K'$. The
   bijection preserves dimensions of underlying simple modules.
  \item[\rm(b)] If $t$ is the restriction to $Y$ of a linear combination
   $\sum_{\chi\in\Irr(KY)} s_\chi \chi$ with $s_\chi\in R_\fp$,
   then the $R'$-linear map $t':Y'\to R' $ is the restriction to $Y'$ of 
   $\sum_{\chi \in \Irr(KY)} \psi(s_\chi)\chi'.$
   Here $\fp$ is the kernel of $\psi$, $R_\fp$ is the corresponding
   localisation and $\psi(s_\chi)$ is the image of $s_\chi$ under
   the unique extension of $\psi$ to a ring homomorphism $R_\fp\to K'$.
 \end{enumerate}
\end{thm}

\begin{proof}
For part (a) see \cite[Thm~68.17, Cor.~68.20] {CR81}. Note that the bijection
given in \cite{CR81} is between the sets $\Irr(\bar K Y)$ and $\Irr(\bar K'Y')$
where $\bar K$ and $\bar K'$ are algebraic closures of $K$ and $K'$
respectively. Since $KY$ and $K'Y'$ are split, this descends via restriction
on both sides to a bijection $\Irr(KY)\to\Irr(K'Y')$. Now~(b) is
an immediate consequence of~(a).
\end{proof}

\section{Hecke algebras and their Schur elements}   \label{sec:Hecke} 
Let $W$ be a finite complex reflection group, that is, a finite subgroup of
$\GL_n(\CC)$ for some $n\ge1$ generated by complex reflections. We denote by
$\QQ_W$ the field generated by the traces of elements of $W$, a finite extension
of $\QQ$. Attached to any reflection $r\in W$ is its reflecting hyperplane
$H=\ker(1-r)$ in $\CC^n$; its point-wise stabiliser $W_H$ in $W$ is cyclic,
generated by reflections. We say that $r\in W_H$ is the \emph{distinguished
reflection} associated to $H$ if $r$ generates $W_H$ and has non-trivial
eigenvalue $\exp(2\pi\bi/o(r))$ of smallest possible argument.

\subsection{From braid groups to trace forms}   \label{subsec:braid}
Let $B(W)$ be the topological braid group of $W$ (see Brou\'e--Malle--Rouquier
\cite{BMR}), that is, the fundamental group of the space of regular orbits of
$W$ on $\CC^n$. So there is an associated exact sequence
$$1\to P(W)\to B(W)\to W\to1,$$
with kernel the pure braid group $P(W)$. Let $A:=\ZZ_W[\bu^{\pm1}]$, where
$\ZZ_W$ denotes the ring of integers of $\QQ_W$ and where $\bu=(u_{rj})$ are
algebraically independent elements indexed by conjugacy classes of distinguished
reflections $r\in W$ and $1\le j\le o(r)$. By \cite[Def.~4.21]{BMR} the
\emph{(generic) Hecke algebra} $\cH(W,\bu)$ of $W$ is defined to be the
quotient of the group algebra $A[B(W)]$ of $B(W)$ over $A$ by the ideal
generated by the elements (called \emph{deformed order relations})
$$ \prod_{j=1}^{o(r)}(\br-u_{rj})\eqno(\cH) $$
for $\br$ running over the braid reflections of $B(W)$ (introduced as
\emph{generators of the monodromy around a hyperplane} in \cite[2B]{BMR}), with
$r$ denoting the image of $\br$ in $W$, a distinguished reflection. We will
write $A[B(W)]\rightarrow\cH(W,\bu)$, $\bx\mapsto h_\bx$, for the associated
quotient map.

We have the following result, first conjectured in \cite{BM93}, with the final
cases having been established by Chavli, Marin and Tsuchioka, respectively, see
\cite{Ch18,Ma19,Ts20} and the references therein:

\begin{thm}[`Freeness Conjecture']   \label{thm:Hecke free}
 The algebra $\cH(W,\bu)$ is $A$-free of rank $|W|$.
\end{thm}

By results of the second author \cite[Cor.~4.8]{Ma99} there is a positive
integer $z$ such that 
the field of fractions $K_W$ of $\tilde A:=\ZZ_W[\tilde\bu^{\pm1}]\supseteq A$
is a splitting field of $\cH(W,\bu)$, where $\tilde\bu=(\tu_{rj})$ are such
that $\tu_{rj}^z=\exp(-2\pi\bi j/o(r))u_{rj}$ for all $r,j$.

Furthermore, by \cite[Prop.~7.1]{Ma00} to each irreducible character $\chi$ of
$K_W\otimes_A\cH(W,\bu)$ is associated an element $f_\chi \in \tilde A$ called
the \emph{Schur element} of~$\chi$. The Schur element
$p_W:=p_W(\bu):=f_{1_W}$ of the trivial character~$1_W$ is called the
\emph{Poincar\'e polynomial} of $\cH(W,\bu)$ (a Laurent polynomial in
the~$\tu_{rj}$). For simply-laced Coxeter groups, this is in fact the
homogenisation of the usual Poincar\'e polynomial of $W$.

The collection $\{(\chi,f_\chi)\mid\chi\in\Irr(\cH(W,\bu))\}$ is
Galois-invariant, so in particular
$$t_{W,\bu}(h) :=\sum_{\chi\in\Irr(\cH(W,\bu))}\frac{1}{f_\chi}\chi(h)$$
lies in $\Frac(A)$ for all $h\in\cH(W,\bu)$. Thus, this defines a symmetric
$A$-linear map
$$t_{W,\bu}: \cH(W, \bu) \to \Frac(A),\quad h \mapsto t_{W, \bu}(h),$$
called the \emph{canonical trace form on $\cH(W,\bu)$}. We denote also by
$t_{W,\bu}$ its $\Frac(A)$-linear extension to
$\cH_{\Frac(A)}(W,\bu):=\Frac(A)\otimes_A\cH(W,\bu)$.
This form satisfies the following property:

\begin{prop}   \label{prop:nice form}
 There exists a set $\cB\subset\cH(W,\bu)$ consisting of monomials in images of
 braid reflections with $1\in\cB$, such that $t_{W,\bu}(b)=\delta_{1,b}$ for
 $b\in\cB$, and $\Frac(A)\otimes\cB$ is a basis of $\cH_{\Frac(A)}(W,\bu)$.
\end{prop}

\begin{proof}
First note that the claim easily reduces to the case of irreducible reflection
groups. For those, it holds by the construction of the Schur elements $f_\chi$,
see \cite{Ma97} and \cite{Ma00} for the exceptional types, and
\cite[Thm~1.3, Lemma~4.3 and \S4.5]{GIM00} for the infinite series.
\end{proof}

Furthermore, $t_{W,\bu}$ satisfies a duality with respect to a certain central
element, but this will not be of importance here. In analogy with the case of
finite Coxeter groups, the theory of spetses predicts the following, see
\cite[Thm-Ass.~2.1]{BMM99}:

\begin{conj}   \label{conj: symmform Hecke}
 The form $t_{W,\bu}$ takes values in $A$ and is a symmetrising form on
 $\cH(W,\bu)$.
\end{conj}

The above conjecture has been established for all imprimitive complex
reflection groups $G(e,1,m)$, for example, by Malle--Mathas \cite[Thm]{MM98}.

Let $W_0\le W$ be a reflection subgroup. We denote by $\cH(W_0,\bu_0)$ the
Hecke algebra of $W_0$ whose parameters $\bu_0$ consist of those parameters
for~ $W$
whose corresponding reflections are, up to conjugacy, contained in $W_0$. Since
non-conjugate reflections in $W_0$ might be conjugate in $W$, $\cH(W_0,\bu_0)$
is a specialisation of the generic Hecke algebra of~$W_0$ corresponding to an
identification of certain of its parameters. It follows from
the Freeness Conjecture (Theorem~\ref{thm:Hecke free}) that $\cH(W_0,\bu_0)$ is
naturally a subalgebra of $\cH(W,\bu)$. Moreover, by the explicit results in
\cite[Thm~5.2]{Ma99} the field $K_W$ is also a splitting field for
$\cH(W_0,\bu_0)$. By $t_{W_0, \bu_0}$ we will mean the corresponding
specialisation of the canonical form of the generic Hecke algebra of~$W_0$. 
Recall that, by a theorem of Steinberg, all parabolic subgroups of $W$, that
is, stabilisers in $W\le\GL_n(\CC)$ of subspaces of $\CC^n$, are reflection
subgroups of $W$.

\begin{defn}   \label{def:strongly quasi}
 We will say that $\cH(W,\bu)$ is \emph{strongly symmetric} if the following
 holds: 
 \begin{enumerate}
  \item[\rm(1)] $t_{W,\bu}$ is a symmetrising form on $\cH(W, \bu)$ and there
   is a section $W\to\bW\subset B(W)$ of the natural map
   $B(W)\to W$ containing~$1$ whose image in $\cH(W,\bu)$ is an $A$-basis
   of $\cH(W,\bu)$ with $t_{W,\bu}(h_\bw) =\delta_{\bw,1}$ for all $\bw\in\bW$;
   and
  \item [\rm(2)] for any parabolic subgroup $W_0\le W$, $t_{W_0,\bu_0}$ is a
   symmetrising form on $\cH(W_0,\bu_0)$ and 
   $t_{W,\bu}|_{\cH(W_0,\bu_0)}=t_{W_0, \bu_0}$.
\end{enumerate}
\end{defn}

Note that the assertion that (2) holds is referred to as the parabolic trace
conjecture in \cite{CC20}. Strong symmetry is known to be satisfied in many
cases; here we use the Shephard--Todd notation for irreducible complex
reflection groups:

\begin{prop}   \label{p:cases}
 For the following irreducible groups, $\cH(W,\bu)$ is strongly symmetric:
 \begin{enumerate}
  \item[\rm(a)] for $W$ a Coxeter group;
  \item[\rm(b)] for $W=G(e,p,n)$ with $n\ne2$ or $p$ odd; and
  \item[\rm(c)] for $W=G_i$, $i\in\{4,5,6,7,8\}$.
 \end{enumerate}
\end{prop}

\begin{proof}
First note that property~(2) of being strongly symmetric follows for a
parabolic subgroup $W_0$ of $W$ if there exists $\bW \subset B(W)$ as in~(1)
such that
\begin{center}
(2') $\{h_\bw\mid w\in W_0\}$ is an $A$-basis of $\cH(W_0,\bu_0)$
and $t_{W_0,\bu_0}(h_\bw) =\delta_{\bw,1}$ for all $w\in W_0$.
\end{center}
For Coxeter groups, (1) and (2') hold for any section $\bW$ consisting of
reduced expressions in the standard generators, see \cite[Prop.~8.1.1]{GP}. For
$W=G(e,p,n)$ the existence of $\bW$ satisfying~(1) is shown in
\cite[Thm~5.1]{MM98}. Since any parabolic subgroup of $W$ is a Young subgroup,
that is, a product of symmetric groups with a group $G(e,p,n')$ for $n'\le n$,
the Ariki--Koike basis of $\cH(W,\bu)$ considered in \cite{MM98} also satisfies
(2'). (See also \cite[p.~177]{BMM99}). The claim for $G_i$,
$i\in\{4,5,6,7,8\}$, follows from the explicit results in \cite{BCCK}.
\end{proof}

\subsection{Specialisations}   \label{subsec:spec}
By a \emph {specialisation} we will mean a ring homomorphism $\psi:A' \to R$
where $A'$ and $R$ are commutative rings with $A'\supseteq A$. We then set
$\cH_\psi(W,\bu):=R\otimes (A'\otimes_A\cH(W,\bu))$. If $\psi$ is inclusion we
will sometimes write
$\cH_R(W,\bu)$ instead of $\cH_\psi(W,\bu)$. The restriction of any
specialisation $\psi$ to a subring of $A'$ will again be denoted by $\psi$ as
will the composition of $\psi$ with any inclusion $R \hookrightarrow R'$.

We will consider certain types of specialisations. For the remainder of this
section, let $R$ be an integral domain containing $\ZZ_W$ and let $K$ be its
field of fractions. Let $\psi_1:R[\tilde\bu] \to R$ be the $R$-linear
homomorphism defined by $\psi_1(\tu_{rj}):=1$ for all $r$ and $j$. So, $\psi_1$
restricts to the specialisation 
$$ A \to \ZZ_W, \quad u_{rj}\mapsto \exp(2\pi\bi j/o(r)). $$
By Bessis \cite[Thm~0.1]{Be01}, the Hecke algebra maps to the group algebra
$RW$ of~$W$ under $\psi_1$. Combining this with the Freeness Conjecture, 
the Tits deformation theorem gives that 
$\cH_{\tilde K} (W,\bu) := \tilde K \otimes_A \cH (W, \bu)$ is isomorphic
to the group algebra $\tilde KW$, where $\tilde K=\Frac(R [\tilde \bu])$
(see Theorem~\ref{t:Tits} and note that $K_W \subseteq \tilde K$).
Thus, we may and will identify the irreducible characters of
$\cH_{\tilde K}(W,\bu) $ with $\Irr(W)$ via the bijection
$\chi_\phi\leftrightarrow\phi$
induced by $\psi_1$. This also induces a labelling of Schur elements of
$\cH(W,\bu)$ by $\Irr(W)$ and we will henceforth denote them as $f_\phi$,
$\phi\in\Irr(W)$. Since $K_W$ is also a splitting field for the Hecke algebra
$\cH(W_0,\bu_0)$ of any reflection subgroup $W_0$ of $W$ we also have
$\cH_{\tilde K}(W_0,\bu)\cong \tilde KW_0$. We will similarly identify the
irreducible characters of $\cH(W_0,\bu_0)$ over $\tilde K$ with $\Irr(W_0)$.

Now let $q$ be a prime power and suppose that $R$ contains $q^{\pm1/z}$. 
We also consider $R$-linear specialisations of the form
$$\psi_q:R[\tilde\bu^{\pm1}]\to R,\quad \tilde u_{rj}\mapsto q^{a_{rj}/z},$$
for integers $a_{rj}$. Any such $\psi_q$ restricts to a specialisation
$$A \to\ZZ_W,\quad u_{rj}\mapsto\zeta_{o(r)}^jq^{a_{rj}},$$
with $\zeta_{o(r)}\in\ZZ_W$ an $o(r)$th primitive root of unity.

\begin{lem}   \label{lem:non-vanishing schur}
 For all $\phi\in\Irr(W)$ we have $\psi_1(f_\phi)=\phi(1)/|W|$ and
 $\psi_q(f_\phi)\ne0$.
 \end{lem}

\begin{proof}
The assertion on $\psi_1$ follows since by construction the set~$\cB$ from
Proposition~\ref{prop:nice form} specialises under $\psi_1$ to the basis of
group elements of $\ZZ_W[W]$. Next, the following can be observed from the
explicit form of the Schur elements (and was first stated explicitly in
\cite[Thm~4.2.5]{ChBook}):
any $f_\phi$ is a product of a scalar, a monomial in the $\tilde u_{rj}$ and
a product of cyclotomic polynomials $\Psi_i$ over $K_W$ evaluated at monomials
$M_i$ in the $\tilde u_{rj}^{\pm1}$ of total degree~0. Thus we need to see that
$\psi_q(\Psi_i(M_i))\ne0$. This is clear if $\psi_q(M_i)$ is not a root of
unity. Now $\psi_q(M_i)$ can only be a root of unity, if the powers of $q$
cancel completely in $\psi_q(M_i)$, which means that $\psi_q(M_i)=\psi_1(M_i)$
and so $\psi_q(\Psi_i(M_i))=\psi_1(\Psi_i(M_i))$. But the latter is a factor of
$\psi_1(f_\phi)=\phi(1)/|W|$ and hence non-zero.
\end{proof}

For the next result we note that symmetrising forms remain symmetrising
after specialisation, that is if $\th:\cO \to \cO'$ is a ring homomorphism,
$Y$ is an $\cO$-algebra which is finitely generated and free as $\cO$-module
and $\tau:Y\to \cO$ is a symmetrising form, then the induced $\cO'$-linear form
$\tau'$ on $Y':= \cO'\otimes_\cO Y$ satisfying $\tau'( 1\otimes y) = \tau(y)$
for $y\in Y$ is a symmetrising form on $Y'$.

\begin{lem}   \label{lem:transfer symm} 
 Assume that $\cH(W,\bu)$ is symmetric over $A$ with respect to the form
 $t_{W,\bu}$ with Schur elements $f_\phi\in\tilde A$, $\phi\in\Irr(\cH(W,\bu))$,
 and let $K=\Frac(R)$. The algebra $K \cH_{\psi_q}(W,\bu)$ is split semisimple.
 Let $t'$ be the induced form on $\cH_{\psi_q}(W,\bu)$. For each
 $\phi\in\Irr(W)$, $\psi_q(f_\phi)$ is the corresponding Schur element of $t'$.
\end{lem}

\begin{proof}
By Lemma~\ref{lem:non-vanishing schur} we have $\psi_q(f_\phi)\ne0$ for all
$\phi\in\Irr(W)$. Now $\psi_q$ is a concatenation of specialisation maps
whose kernel is a prime ideal of height~1. So by \cite[Thm~2.4.12]{ChBook},
$K\cH_{\psi_q}(W,\bu)$ is also split semisimple. The result thus follows by
Tits' deformation theorem (Theorem~\ref{t:Tits}).
\end{proof}

Finally, for an indeterminate $x$ we will also consider the \emph{spetsial
specialisation}
$$\psi_\spets: R[\tilde\bu^{\pm1}] \to R[x^{\pm\frac{1}{z}}], \qquad
  \tilde u_{rj}\mapsto\begin{cases} x^{\frac{1}{z}}& \text{if $j=o(r)$,}\\
   1& \text{if $1\le j<o(r)$}\end{cases}$$
(so $\psi_\spets(u_{r,o(r)})=x$, $\psi_\spets(u_{rj})=\exp(2\pi\bi j/o(r))$
for $j<o(r)$). Clearly $\psi_1$ factorises through $\psi_\spets$ by composition
with $x^{1/z}\mapsto1$ as does any specialisation $\psi_q$ for the case
$a_{rj} = 1$ if $j =o(r)$ and $a_{rj} = 0$ if $j< o(r)$ by composition with
$x^{1/z}\mapsto q^{1/z}$. We write $\psi_{\spets,q}$ for this latter
specialisation, which will become important in Section~\ref{subsec:pf thm1}.
The spetsial specialisation links Schur elements of Hecke algebras to
unipotent character degrees of spetses (see Section~\ref{subsec:spetsial}).

\section{Conjectures for the principal block of a $\ZZ_\ell$-spets}   \label{sec:sec3}
In this section we define the dimension of the principal block of a
$\ZZ_\ell$-spets, introduced in \cite[\S6.2]{KMS}, and propose some
conjectures around this notion.

\subsection{The principal block of a $\ZZ_\ell$-spets} \label{subsec:B_0}
Let $\ell $ be a prime and let $q$ be a prime power not divisible by~$\ell$.
Recall that (under some conditions) the set of
characters of the principal $\ell$-block of a finite reductive group $G$ over
$\FF_q $ can be described as a union of sets of characters in bijection with
the principal $e$-Harish Chandra series of unipotent characters of
dual-centraliser subgroups $C_{G^*}(s)^*$ where $G^*$ is the ``Langlands dual"
of $G$ and $s$ runs over conjugacy classes of $\ell$-elements of $G^*$. The
principal block of a $\ZZ_\ell$-spets as constructed in \cite{KMS} is modelled
on this description: unipotent characters and the appropriate Harish-Chandra
series are provided by the theory of spetses whereas the indexing set of
$\ell$-elements and corresponding centralisers comes from the theory of
$\ell$-compact groups and fusion systems. In the next paragraph, we briefly
recall this construction. Before doing so, we point out that the description of
the principal block of a $\ZZ_\ell$-spets does not involve going over to the
dual as one would expect from analogy with the group case. The reason for this
is that the fusion system construction that we rely on is for the moment only
available for simply connected $\ell$-compact groups. However, this departure
does not lead to an inconsistency in the group case under the conditions on
$\ell$ with which we are concerned (see \cite[Prop.~6.8]{KMS}). 

We assume from now till the end of the section that $\ell>2$. Let $W\le\GL(L)$
be a finite \emph{spetsial} (see
\cite[\S3]{MaICM}) $\ell$-adic reflection group on a $\ZZ_\ell$-lattice $L$,
and let $X$ be the associated connected $\ell$-compact group (see
\cite[Thm 1.1]{BM07}). We assume moreover that $X$ is simply connected and
that $\ell$ is \emph{very good} for $(W,L)$, in the sense of
\cite[Def.~2.4]{KMS}. Let $\vhi\in N_{\GL(L)}(W)$ be of $\ell'$-order and
$\GG=(W\vhi^{-1},L)$ be the associated simply connected $\ZZ_\ell$-spets. For
example any Weyl group $W$ for which $\ell$ is very good (in the classical
sense) determines a $\ZZ_\ell$-spets satisfying the above conditions.

Set $\GG(q):=(W\vhi^{-1},L,q)$. In \cite{BM07}, Broto and
M\o ller showed how to attach to these data a fusion system $\cF$ on a finite
$\ell$-group $S$, via the associated $\ell$-compact group $X$
(see~\cite[Thm~3.2]{KMS}). Here, $S$ is an extension of a homocyclic
$\ell$-group $T$ (the \emph{toral part}) by a Sylow $\ell$-subgroup of the
associated relative Weyl group $W_{\vhi^{-1}\zeta^{-1}}$ (see
\cite[Thm~3.6]{KMS}). Note that $S$ and $\cF$ only
depend on the $\ell$-part $\ell^a$ of $q-\zeta$, where
$\zeta\in\ZZ_\ell^\times$ is the root of unity with $q\equiv\zeta\pmod\ell$,
not on $q$ itself. 

If $\GG(q)$ arises from a connected reductive algebraic group $\bG$ over
$\overline{\FF_q}$ with Weyl group $W$ and a Frobenius morphism $F:\bG\to\bG$
with respect to an $\FF_q$-structure acting as $\vhi$ on $W$, then $S$
is a Sylow $\ell$-subgroup of $\bG^F$ and $\cF$ is the $\ell$-fusion system of
$\bG^F$ on $S$ (see \cite[Rem.~3.3 and Sec.~5.3]{KMS}). In particular,
$T=(\bT^F)_\ell$ where $\bT$ is a maximal $e$-split torus of $\bG$.

Under our assumptions, for any $s\in S$ the centraliser $W(s):=C_W(s)$ of $s$
is again an $\ell$-adic reflection group, a reflection subgroup of $W$ (see the
proof of \cite[Thm~5.2]{KMS} or \cite[Prop.~2.3]{KMS})), and by
\cite[Prop.~6.2]{KMS}, it is again spetsial. We let
$C_\GG(s):=(W(s)\vhi_s^{-1},L)$ be the associated $\ZZ_\ell$-spets, where
$\vhi_s\in N_{W\vhi}(W(s))$ is defined as in \cite[\S5.2]{KMS}.

Now recall from \cite{Ma95} (for the infinite series of irreducible complex
reflection
groups) and \cite{BMM14} (for the primitive ones) that associated to $\GG$ as
well as to the various $C_\GG(s)$ there are sets of unipotent characters
$\Uch(\GG)$ and $\Uch(C_\GG(s))$, respectively. If $W$ is a Weyl group, these
are just the unipotent characters of an associated finite reductive group. For
$s\in S$ we let
$$\cE(\GG,s)=\{\gamma_{s,\la}\mid\la\in\Uch(C_\GG(s))\}$$
denote a set in bijection with $\Uch(C_\GG(s))$ and call it the
\emph{characters of $\GG$ in the series~$s$}. The sets $\Uch(C_\GG(s))$ are
in canonical bijection for conjugate elements $s$. Moreover, these sets only
depend on $\ell^a$, not on $q$ itself. Any unipotent character $\la$ comes with
a degree (polynomial) $\la(1)\in\ZZ_\ell[x]$. The \emph{degree of
$\gamma_{s,\la}\in\cE(\GG,s)$} is defined as
$$\gamma_{s,\lambda}(1):=|\GG:C_\GG(s)|_{x'}\,\la(1)\in\ZZ_\ell[x].$$
Here, $|\GG:C_\GG(s)|_{x'}$ means the prime-to-$x$ part of the polynomial
$|\GG|/|C_\GG(s)|\in\ZZ_\ell[x]$, where $|\GG|,|C_\GG(s)|$ are the respective
order polynomials; note that the latter divides the former by
\cite[Lemma~6.6]{KMS}.

Now by \cite[Folg.~3.16 and~6.11]{Ma95} and \cite[4.31]{BMM14} for any
$\ZZ_\ell$-spets $\HH$, for any root of unity $\eta$ the set of unipotent
 characters $\Uch(\HH)$ is naturally partitioned into so-called
$\eta$-Harish-Chandra series, and one among them, the \emph{principal
$\eta$-Harish-Chandra series $\cE(\HH,1,\eta)$} of $\Uch(\HH)$ containing
$1_\GG$, is in bijection with the irreducible characters of the corresponding
Springer--Lehrer relative Weyl group. In particular, $\cE(\HH,1,1)$ is in
bijection with $\Irr(W^\vhi)$.

With this, for $\zeta$ as above denote by $\cE(\GG,s)_\zeta$ the subset of
$\cE(\GG,s)$ in bijection with the principal $\zeta$-Harish-Chandra series of
$\Uch(C_\GG(s))$, and hence also in bijection with the irreducible characters of
the relative Weyl group $W(s)_{\vhi_s^{-1}\zeta^{-1}}$.

The following definition from \cite[\S6.2]{KMS} is inspired by the results of
Cabanes--Enguehard on unipotent $\ell$-blocks of finite reductive groups;
indeed, if $\GG$ is a rational spets for which $\ell$ is very good, then what
we define are exactly the characters in the principal $\ell$-block of the
corresponding finite group of Lie type $\GG(q)$ (see \cite[Prop.~6.8]{KMS}):

\begin{defn}   \label{def:B0}
 Let $\GG=(W\vhi^{-1},L)$ be a simply connected $\ZZ_\ell$-spets with $\vhi$ of
 $\ell'$-order such that $\ell$ is very good for $\GG$ (in the sense of
 \cite[Def.~2.4]{KMS}) and $q$ a prime power with $q\equiv\ze\pmod\ell$. The
 \emph{characters in the principal block $B_0$ of $\GG(q)$} are
 $$\Irr(B_0):=\coprod_{s\in S/\cF}\cE(\GG,s)_\zeta,$$
 where the union runs over a set $S/\cF$ of representatives $s$ of
 $\cF$-conjugacy classes in $S$. The \emph{dimension of $B_0$} is defined as
 $$\dim(B_0)=\sum_{\gamma\in\Irr(B_0)}\gamma(1)^2
   =\sum_{s\in S/\cF}\sum_{\la\in\cE(C_\GG(s),1,\zeta)} \gamma_{s,\la}(1)^2
   \in\ZZ_\ell[x].$$
\end{defn}

Again, these do not depend on $q$ but only on $\ell^a=(q-\zeta)_\ell$. 

The dimension of the principal $\ell$-block $B_0$ of a finite group $G$ with
Sylow $\ell$-subgroup $S$ satisfies $(\dim B_0)_\ell=|S|$ (see 
Proposition~\ref{prop:dimgreen}). We conjecture that our dimension of the
principal block $\dim(B_0)$ of $\GG(q)$ behaves similarly:

\begin{conj}   \label{conj:dim B}
 Let $\GG=(W\vhi^{-1},L)$ be a simply connected $\ZZ_\ell$-spets such that
 $\ell$ is very good for $\GG$. Let $q$ be a prime power not divisible by
 $\ell$ and $S$ the associated $\ell$-group. Then
 $$(\dim(B_0)|_{x=q})_\ell = |S|.$$
\end{conj}

By further analogy with the group case (see Proposition~\ref{prop:dimgreen})
we also conjecture the following global-local statement:

\begin{conj}   \label{conj:quot}
 In the setting of Conjecture~{\rm\ref{conj:dim B}}, if
 $\zeta\in\ZZ_\ell^\times$ with $q\equiv\zeta\pmod\ell$ we have
 $$\big(\dim(B_0)|_{x=q})_{\ell'}
   \equiv |W_{\vhi^{-1}\ze^{-1}}|_{\ell'}\pmod\ell\,.$$
\end{conj}

Note that Conjectures \ref{conj:dim B} and~\ref{conj:quot} combine to form
Conjecture \ref{c:main}.

\begin{exmp}   \label{ex:specialspets}
We describe $S$, $\cF$, and $\Irr(B_0)$ in a special case relevant to
Theorem~\ref{t:main}. Suppose that $q\equiv 1\pmod\ell$ and $\varphi=1$. Then
$T$ may be identified with $L/\ell^a L$ where $\ell^a ||(q-1)$,
$W = W_{\vhi^{-1}\zeta^{-1}}$ and the action of $W$ on $T$ in $S$ is the one
inherited from the action of $W$ on $L$. Assume in addition that $W$ is an
$\ell'$-group. Then the condition that $\ell$ is very good for $(W,L)$ always
holds (see \cite[Prop.~2.6]{KMS}). Further, $S=T$ and $\cF$ is the
$\ell$-fusion system $\cF_{SW}(S)$ of the group $SW$ on $S$ (see
\cite[Thm~3.4]{KMS} and \cite[Thm~9.8]{BM07}). Moreover, for any $s \in S=T$
the subgroup $W(s)_{\vhi_s^{-1}\zeta^{-1}}$ equals $W(s)$ and hence
$\cE(\GG,s)_1$ is in bijection with $\Irr(W(s))$. So $\Irr(B_0)$ is in
bijection with $W$-classes of pairs $(s,\phi)$ where $s\in S$ and
$\phi\in\Irr (W(s))$. 
\end{exmp}

\subsection{Decomposition numbers}   \label{subsec:deco}
Suppose in this subsection that $|W|$ is prime to $\ell$, that $\ell\mid(q-1)$
and that $\varphi=1$. Recall the description of the principal block $B_0$ under
these assumptions given in Example~\ref{ex:specialspets}.

Since $W$ is an $\ell'$-group and $S$ an $\ell$-group, we may identify
$\IBr(SW)$ with the subset $\Irr(W)$ of $\Irr(SW)$. Similarly, we think
of the unipotent characters in $B_0$ as the irreducible Brauer characters of
$B_0$ and set
$$ \IBr(B_0) :=\cE(\GG, 1)_1 \subseteq \Irr(B_0). $$
We associate decomposition numbers, and formal degrees of projective
indecomposable characters to $B_0$ as follows.

The $\cF$-classes of elements of $S$ are the $W$-conjugacy classes of $S$.
Further,
since $(|W|,\ell) = 1$, the Glauberman--Isaacs correspondence gives that the
actions of $W$ on $S$ and on $\Irr(S)$ are permutation isomorphic. Thus there
is a bijection between the set of $W$-classes of $\Irr(S) $ and the set of
$W$-classes of $S$ such that if the class of $s\in S$ corresponds to the
class of $\hat s\in\Irr(S)$, then $ W_s= W_{\hat s}$ where $W_s, W_{\hat s}$
denotes the stabiliser in $W$ of $s$, $\hat s$ respectively. Note that $W_s$
was denoted $W(s)$ in Example~\ref{ex:specialspets}.
Such a bijection between the set of $W$-classes of $\Irr(S)$
and of $S$ will be called \emph{$W$-equivariant} if in addition the (class of)
$1\in S$ is sent to the (class of the) trivial character of $S$. Note that a
$W$-equivariant bijection always exists.

By Clifford theory,
$$ \Irr(SW)= \coprod_{\th \in \Irr(S)/W}\Irr(SW|\th),$$
where the union runs over a set $\Irr(S)/W$ of representatives $\th$ of
$W$-conjugacy classes of $\Irr(S)$ and where $\Irr(SW|\th)$ denotes the
set of irreducible characters of $SW$ covering~$\th$. Moreover, since $|W|$
is prime to $\ell$, $\Irr(SW|\th)$ is in
bijection with $\Irr(W_\th)$.

By the description of $\Irr(B_0)$ given in Example~\ref{ex:specialspets},
$|\Irr(B_0)| =|\Irr(SW)|$. A bijection
$\Theta: \Irr(SW)\to\Irr(B_0)$, $\gamma\mapsto\hat\gamma$, will be said to be
$W$-equivariant if there exists a $W$-equivariant bijection
$\Irr(S)\to S$ such that for corresponding elements $s\in S$ and
$\hat s\in\Irr(S)$, $\Theta$ restricts to a bijection
$\Irr(SW|\hat s)\to\cE(\GG, s)_1$. Since $\IBr(SW)$ is in bijection with
$\Irr(W|1)$, any $W$-equivariant bijection $\Irr(SW)\to\Irr(B_0)$
restricts to a bijection $\IBr(SW)\to\IBr(B_0)$.
 
Let $\Irr(SW)\to\Irr(B_0)$, $\gamma\mapsto\hat\gamma$, be a $W$-equivariant
bijection. We
declare the decomposition matrix of $B_0$ to be the $\ell$-decomposition matrix
of $SW$ via this bijection, that is if $d_{\gamma\nu}$ is the decomposition
number in $SW$ corresponding to $\gamma\in\Irr(SW)$ and $\nu\in\IBr(SW)$,
then we regard $d_{\gamma\nu}$ also as the decomposition number for
$\hat\gamma\in\Irr(B_0)$ and $\hat\nu\in\IBr(B_0)$. Recall that for any
$\gamma\in\Irr(SW)$, we have
$\gamma(1) = \sum_{\nu\in \IBr(SW)}d_{\gamma\nu}\nu(1)$. In
Proposition~\ref{prop:Brauerdegreeconsistent} we show that the analogous
equations hold in $B_0$. We define 
$$\deg\Phi_{\hat \nu}
  :=\sum_{\gamma \in\Irr(SW)}d_{\gamma \nu}\deg(\hat \gamma)\in\ZZ_\ell[x]
  \qquad\text{for $\hat \nu\in\IBr(B_0)$},$$
to be the \emph{formal degrees of projective indecomposable characters of
$B_0$}. The following is a restatement of Conjecture \ref{c:main2}.

\begin{conj}   \label{conj:decomp}
 In the setting of Conjecture~{\rm\ref{conj:dim B}}, if $W$ has order coprime
 to $\ell$, $\varphi $ is trivial and $q\equiv1\pmod\ell$, then for some
 $W$-equivariant bijection $\Irr(SW)\buildrel\sim\over\to\Irr(B_0)$ we have
 that $|S|$ divides
 $(\deg\Phi_{\hat \nu})|_{x =q}$ for all $\hat\nu\in\IBr(B_0)$.
\end{conj}

\subsection{The rational and the primitive cases}
We'll prove the above conjectures for most $W$ of order coprime to $\ell$ in
Theorem~\ref{thm:omni}.
For the moment, let us see why they hold in the rational case:


\begin{prop}   \label{prop:rational}
 Conjectures~{\rm\ref{conj:dim B}}, {\rm\ref{conj:quot}}, and
 {\rm\ref{conj:decomp}} hold if $\GG$ is a rational spets underlying a
 finite reductive group.
\end{prop}

\begin{proof}
Let $\bG$ be a connected reductive group over an algebraically closed field of
characteristic~$p$ and $F:\bG\to\bG$ a Frobenius endomorphism with respect to
an $\FF_q$-structure, such that $\GG$ is the underlying spets. That is, $\vhi$
is the automorphism of $W$ induced by~$F$. Recall that $S$ may be identified
with a Sylow $\ell$-subgroup of $\bG^F$ and $\cF$ with the fusion system
$\cF_{\bG^F}(S)$ (see \cite[Rem.~3.3(a) and Sec.~5.3]{KMS}). Let $d$ be the
order of $\zeta$, hence the order of $q$ modulo~$\ell$. By
\cite[Prop.~6.8]{KMS}, there is a degree preserving bijection between
$\Irr(B_0)$ and the set of irreducible characters of the principal $\ell$-block
$B_0(\bG^F)$.
Note that the bijection given in \cite[Prop.~6.8]{KMS} is stated to preserve
defects but it is easy to check from the setup that for any irreducible
character $\gamma_{s,\la}$ in $B_0$, $\gamma_{s,\la}(1)|_{x=q}$ is the degree
of the corresponding character of $B_0(\bG^F)$. Now the assertion regarding
Conjecture~\ref{conj:dim B} follows from Proposition~\ref{prop:dimgreen}.

As described in Section~\ref{subsec:B_0}, $S = T.(W_1)_\ell$ with
$W_1:=W_{\vhi^{-1}\ze^{-1}}$. Under our assumptions on~$\ell$, $\bL:=C_\bG(T)$
is a Levi subgroup, and moreover $N_\bG(S)^F\leq N_\bG(T)^F=N_\bG(\bL)^F$ (see
\cite[Thm~5.9 and~5.14]{Ma07}). Let $H:=N_\bG(\bL)^F$. Then 
Proposition~\ref{prop:dimgreen} shows that in order to prove
Conjecture~\ref{conj:quot} for $\GG$ it suffices to see that
$(\dim B_0(H))_{\ell'}\equiv |W_1|_{\ell'}\pmod\ell$. Now $B_0(H)$ is
isomorphic to the principal block of $H/O_{\ell'}(\bL^F)$. Since $S/T$
acts faithfully on $T$, $T$ is a Sylow $\ell$-subgroup of $\bL^F$, so
$H/ O_{\ell'} (\bL^F)\cong T(N_\bG(\bL)^F/\bL^F)\cong T W_1$, the latter
by the definition of relative Weyl groups.
The result follows as $T W_1$ has a unique $\ell$-block.

Finally, we prove Conjecture \ref{conj:decomp} in this situation. So assume
$\varphi=1$ and $q\equiv 1\pmod\ell$. Then $W_1 = W$. As recalled above, there
is a degree preserving bijection $\Irr(B_0)\to\Irr(B_0(\bG^F))$. On the other
hand, by a result of Puig \cite[Thm~5.5, Cor.~5.10]{Pu94}, as explained in
\cite[Prop.~8.11]{Ca18}, the principal block of $\bG^F$ over a suitably large
complete discrete valuation ring $\cO$ of characteristic~$0$, is Morita
equivalent to the group algebra $\cO[SW]$ (note here $S=T$ as $|W|$ is
prime to $\ell$ and that $W = N_\bG(\bT)^F/\bT^F $, where $\bT$ is a
Sylow $1$-torus of $\bG$ with $T$ the Sylow $\ell$-subgroup of $\bT^F$). In
particular the decomposition
numbers of $\cO[SW]$ are decomposition numbers of $B_0(\bG^F)$. Thus
$(\deg\Phi_{\hat \nu})|_{x=q}$ is the dimension of a projective indecomposable
module of $B_0(\bG^F)$. Now Conjecture \ref{conj:decomp} follows since the
dimension of any projective indecomposable module of a finite group algebra
$\cO G$ is divisible by $|G|_\ell$ (for instance, apply
Proposition~\ref{lem:projschur} with respect to the standard symmetrising form
on~$\cO G$.)
\end{proof}

To deal with the primitive cases, we use the following:

\begin{lem}   \label{l:orlik}
 Let $(W,L)$ be a finite $\ell$-adic reflection group with $|W|$ prime to $\ell$
 and let $W_0\le W$ be a parabolic subgroup. For $1 \le k \le \rk(W_0)$ there
 exist $b^{W_0}_k \in \ZZ$ such that for any prime power $q\equiv1 \pmod\ell$
 and $\cF$ the fusion system attached to $(W,L,q)$, on a homocyclic
 $\ell$-group $T$ of exponent $a$, the number of $W$-orbits ($\cF$-classes) of
 elements of $T$ with stabiliser conjugate to $W_0$ is given by
 $$\frac{1}{|N_W(W_0):W_0|} \prod_{k=1}^{\rk(W_0)} (\ell^a-b^{W_0}_k).$$
\end{lem}

\begin{proof}
Let $\cA$ denote the set of 1-eigenspaces of reflections in $T$ and denote by
$\cL=\cL(\cA)$ the lattice of all intersections of elements of $\cA$ with
minimal element $T$. By Steinberg's theorem (see \cite[Prop.~2.3]{KMS}), each
$Y\in \cL$ is the centraliser in $T$ of some parabolic subgroup. Thus by
inclusion/exclusion, the number of $W$-orbits of elements of~$T$ with
stabiliser conjugate to $W_0$ is given by the Euler characteristic of the
sublattice $\{Y\in\cL\mid C_T(W_0)\le Y\}$ divided by $|N_W(W_0):W_0|$. This
Euler characteristic has the stated form by \cite[Thm~1.2]{OS82}. 
\end{proof}

The tables in \cite{OS82} explicitly list the integers $b^{W_0}_k$ (and the
quantities $|N_W(W_0):W_0|$) for all parabolic subgroups $W_0$ of all
exceptional complex reflection groups $W$. An immediate consequence is the
following result:

\begin{prop}\label{prop:prim}
 Conjectures~{\rm\ref{conj:dim B}} and~{\rm\ref{conj:quot}} hold for all
 primitive spetsial $\ell$-adic reflection groups with $q\equiv 1\pmod\ell$.
\end{prop}

\begin{proof}
If $W$ is a Weyl group, the claim follows from Proposition~\ref{prop:rational}.
For the remaining Coxeter groups, it follows from Theorem \ref{thm:omni}.
Otherwise since $\ell$ is very good, we must have $\ell \nmid |W|$ and
$$W \in\{G_4,G_6,G_8,G_{14},G_{24},G_{25},G_{26},G_{27},G_{29},G_{32},
  G_{33},G_{34}\}.$$
For all of these only $\vhi=1$ is possible, by \cite[Prop.~3.13]{BMM99}. We
explicitly calculate $\dim(B_0)$ as a polynomial in $x$ using
Lemma~\ref{l:orlik} and the tables in \cite[App.]{BMM14}. The required
congruences for $\dim(B_0)|_{x=q}$ are readily checked via the substitution
$q\mapsto 1+r\ell^a$ for $r\in\ZZ$.
\end{proof}

\subsection{Character degrees and Schur elements.}   \label{subsec:spetsial}
The proof of Theorem~\ref{t:main} goes through the connection between unipotent
character degrees of spetses and Schur elements of corresponding Hecke
algebras. We describe this connection in the relevant special case.
Let $\GG=(W,L)$ be a simply connected $\ZZ_\ell$-spets such that
$\ell$ is very good for $\GG$ and $q$ a prime power and let
$\psi_\spets: \ZZ_\ell[\tilde\bu^{\pm1}]\to \ZZ_\ell[x^{1/z}]$
be the spetsial specialisation described in Section~\ref{subsec:spec} with
$R=\ZZ_\ell$.
The degrees of the unipotent characters $\Uch(\GG)$ of $\GG$ in the principal
1-Harish-Chandra series $\cE(\GG,1)_1$ are given by
$$\psi_\spets (f_{1_W}) /\psi_\spets (f_\phi) \qquad
 \text{for $\phi\in\Irr(W)$}$$
(see \cite[S\"atze~3.14, 6.10]{Ma95} for the infinite series and
\cite[Ax.~4.16]{BMM14} for the exceptional types). This leads to the following
formula for character degrees in the principal block $B_0$ of $\GG(q)$.

\begin{lem}   \label{lem:degs1}
 Assume $q\equiv 1\pmod\ell$. Then the degrees of the characters in
 $\Irr(B_0)$ are given by
 $$\coprod_{s\in S/\cF}
   \big\{\frac{|C_\GG(\TT)|_{x'}}{|C_{C_\GG(s)}(\TT_s)|_{x'}}
   \frac{\psi_\spets(p_{W})}{ \psi_\spets (f_{s,\phi})}\mid \phi\in\Irr(W(s))\big\} $$
 where $\TT$, $\TT_s$ is a Sylow $1$-torus of $\GG$, $C_\GG(s)$, respectively,
 and the $f_{s,\phi}$ denote the Schur elements of the Hecke algebra
 $\cH(W(s),\bu_s)$ of $W(s)$, with the parameters $\bu_s$ inherited from
 $\cH(W,\bu)$.
\end{lem}

\begin{proof}
As mentioned above for any $s\in S$ the degrees of the unipotent characters in
the principal $1$-series of $\cE(C_\GG(s),1)$ are given by
$$\big\{\psi_\spets(p_{W(s)})/\psi_\spets (f_{s,\phi})\mid
  \phi\in\Irr(W(s))\big\}\subset\ZZ_\ell[x].$$
Now by \cite[Prop.~8.1]{Ma00} we have
$\psi_\spets(p_W)=\psi_\spets(f_1)=|\GG:C_\GG(\TT)|_{x'}$ and accordingly
$\psi_\spets(p_{W(s)})=|C_\GG(s):C_{C_\GG(s)}(\TT_s)|_{x'}$. Since by
definition the degree of $\gamma_{s,\phi}\in\Irr(B_0)$ is
$$|\GG:C_\GG(s)|_{x'}\ \psi_\spets(p_{W(s)})/\psi_\spets (f_{s,\phi}),$$
our claim follows.
\end{proof} 

\begin{lem}   \label{lem:degree schur base}
 In the situation of Lemma~{\rm\ref{lem:degs1}}, assume moreover that $|W|$ is
 coprime to~$\ell$. Then the degrees of the characters in $\Irr(B_0)$ are
 given by
 $$\coprod_{s\in S/\cF}\big\{\psi_\spets(p_W)/\psi_\spets(f_{s,\phi})
   \mid\phi\in\Irr(W(s))\big\},$$
 where $f_{s,\phi}$ denotes the Schur elements of the Hecke algebra
 $\cH(W(s),\bu_s)$. 
  \end{lem}

\begin{proof}
If $\ell$ does not divide $|W|$ then we have $S=T$, that is, the centraliser of
any $\ell$-element $s\in S$ contains the Sylow 1-torus $\TT$, whence
$\TT_s=\TT$ for all $s$. Now the centraliser of a Sylow 1-torus is a maximal
torus since the coset $W\phi$ always contains a 1-regular element by
\cite[Prop.~3.3]{Ma06}; for this note that none of the exceptions in loc.~cit.
is spetsial. So in fact $C_\GG(\TT)=C_{C_\GG(s)}(\TT)$ and the stated formula
follows from Lemma~\ref{lem:degs1}.
\end{proof}

By analogy with the group case, we now establish a Brauer reciprocity formula
for the Brauer characters and decomposition numbers defined in
Section~\ref{subsec:deco}.

\begin{prop}   \label{prop:Brauerdegreeconsistent}
 Suppose that $W$ is $\ell$-adic spetsial of order coprime to $\ell$, $\vhi=1$
 and $q\equiv1\pmod\ell$. Then for any $W$-equivariant bijection $\hat{\ }:
 \Irr(SW)\to\Irr(B_0)$ we have 
 $$\deg(\hat\gamma) = \sum_{\chi\in\IBr(SW)} d_{\gamma,\chi}\deg(\hat\chi). $$
\end{prop}

\begin{proof}
By Clifford theory the ordinary irreducible characters of $SW\!=\!TW$ are
obtained~as
$$\Irr(TW)=\{\Ind_{TW_\th}^{TW}(\th\otimes\nu)\mid
   \theta\in\Irr(T),\ \nu\in\Irr(W_\th)\}.$$
Since $T=O_\ell(TW)$ and $|W|$ is prime to $\ell$, $\IBr(TW)$ consists
of the restrictions to $\ell'$-classes of the irreducible characters
$1\otimes\nu$, $\nu\in\Irr(W)$, and we may (and will) thus identify $\IBr(TW)$
with $\Irr(W)$. 
The $\ell$-decomposition numbers of $TW$ are then described as follows:
If $\eta\in\Irr(TW)$ then its restriction to $\ell'$-classes $\eta^0$ can be
considered as character of $W$, and
the multiplicity of $\chi\in\IBr(TW)=\Irr(W)$ in $\eta^0$ is just
$\langle\eta,\chi\rangle$. That is, if $\eta=\Ind_{TW_\th}^{TW}(\th\otimes\nu)$
as above, then this multiplicity is $\langle\nu,\chi|_{W_\th}\rangle$.
\par
Now assume that $\gamma\in\Irr(B_0)$. Then there is $s\in T$ and
$\la\in\Uch(C_\GG(s))_1$ such that $\gamma=\gamma_{s,\la}$. Let
$\phi\in\Irr(W(s))$ be the irreducible character indexing
$\la\in\Uch(C_\GG(s))_1$. Now by Lemma~\ref{lem:degree schur base} we have
$\gamma_{s,\la}(1)=\psi_\spets(p_W)/\psi_\spets(f_{s,\phi})$. On the other
hand, for $\gamma'\in\cE(\GG,1)_1$ labelled by $\chi\in\Irr(W)$ we have
$\gamma'(1)=\psi_\spets(p_W)/\psi_\spets(f_\chi)$. Thus the required equality
reads
$$\psi_\spets(f_{s,\phi})^{-1}
   =\sum_{\chi\in\Irr(W)}\langle\phi,\chi|_{W_s}\rangle
\psi_\spets(f_\chi)^{-1}.$$
But this holds for spetsial $W$ by the validity of 1-Howlett--Lehrer
theory, see \cite[S\"atze~3.14, 6.10]{Ma95} for the infinite series and
\cite[Ax.~4.16]{BMM14} for the exceptional types.
\end{proof}

\section{Yokonuma type algebras for torus normalisers for $\ell$-adic reflection groups}   \label{sec:Yoku}
When the order of $W$ is prime to $\ell$ our definition of principal block and
Conjectures~\ref{conj:dim B}, \ref{conj:quot} and~\ref{conj:decomp} are related
to a generalisation of Yokonuma algebras to $\ell$-adic reflection groups.
The classical Yokonuma algebra was defined as the endomorphism algebra of the
permutation representation of a finite Chevalley group on a maximal unipotent
subgroup \cite{Yo67}. It is a deformation of the group algebra of the
normaliser of a maximally split torus, to which it becomes isomorphic over a
splitting field (see \cite[\S34]{Lu7}).
We propose to extend this construction over the $\ell$-adic integers to ``torus
normalisers'' arising from $\ell$-compact groups attached to arbitrary
$\ell$-adic reflection groups. This will allow us in
Section~\ref{subsec:pf thm1} to prove Conjectures~\ref{conj:dim B},
\ref{conj:quot} and~\ref{conj:decomp} in
the case $q\equiv1\pmod\ell$ and $\vhi=1$.

\subsection{Definition and first properties}
Let $\ell$ be a prime and $W$ be a finite $\ell$-adic reflection group, that is,
$W\le\GL(L)$ with $L=\ZZ_\ell^n$. Let $q$ be a prime power with
$q\equiv1\pmod\ell$ and $a$ the positive integer such that $\ell^a||(q-1)$.
Let $T=L/\ell^aL$. Then $T$ is homocyclic of exponent~$\ell^a$ and is equipped
with a natural action of~$W$. For any reflection $r\in W$ we set
$T_r:=[T,r]:=\langle [t,r]\mid t\in T\rangle\le T$.

The topological braid group $B:=B(W)$ of $W$ (see Section~\ref{subsec:braid}) 
acts naturally on $T$ through its quotient $W$. We let $\hB$ be the semidirect
product of $T$ with $B$. Observe that $P(W)$ acts trivially on $T$, so $\hB$
is a (non-split) extension of $T\times P(W)$ by $W$.
\par

Recall from Section~\ref{subsec:braid} the indeterminates $\bu=(u_{rj})$
attached to $W$. We define a new set $\bv=(v_{rj})$ of indeterminates by the
linear relations
$$u_{rj}=\zeta_{o(r)}^j(1+|T_r|v_{rj})\qquad\text{for $r\in W,\ 1\le j\le o(r)$,}$$
where, for any $k|(\ell-1)$, $\ze_k\in\ZZ_\ell$ denotes a primitive $k$th root
of unity. Let $\hA=\ZZ_\ell[\bv,\bu^{-1}]$.
\begin{defn}
 Define $\cY(W,a,\bv)$ to be the quotient of the group algebra $\hA[\hB]$ of
 $\hB=T\rtimes B$ over $\hA$ by the ideal generated by the deformed order
 relations
 $$\prod_{j=1}^{o(r)}\big(\br-\ze_{o(r)}^j(1+v_{rj}E_r)\big)\qquad
   \text{with $E_r:=\sum_{t\in T_r}t\in \hA[T]$,}\eqno(\dagger)$$
 where $\br$ runs over the braid reflections of $B\le\hB$ and~$r$ denotes
 the image of $\br$ in $W$. We will write $x\mapsto y_x$ for the canonical map
 $\hat A[\hB]\to\cY(W,a,\bv)$.
\end{defn}

When $W$ is a Weyl group, the deformed order relation $(\dagger)$ generalises
the quadratic one from the classical Yoko\-numa algebra
\cite[Thm~2(2.1)]{JK01}, see also \cite[2.2(3)]{Ma18a}. In
Section~\ref{sub:classical yokonuma}, we show that in this case a suitable
specialisation of $\cY$ is isomorphic to a truncation of the classical Yokonuma
algebra.

As far as we can tell, there is no direct relation between the algebra
$\cY(W,a,\bv)$ defined above and the ``cyclotomic Yokonuma--Hecke algebra''
considered by Chlouveraki--d'Andecy \cite[\S2]{CdA16} for the reflection
group $W=G(d,1,n)$; in their algebra, the underlying reflection group
$W$ only acts via its quotient $G(1,1,n)\cong\fS_n$ on the torus. On the
other hand, our construction is related to an algebra defined by Marin
\cite{Ma18a}, see Remark~\ref{rem:Marin} below.

Henceforth, for simplicity we set $\cY:=\cY(W,a,\bv)$. Note that the
specialisation $\psi_1:\ZZ_\ell[\bu^{\pm1}] \to \ZZ_\ell$,
$u_{rj}\mapsto\ze_{o(r)}^j$, extends to a homomorphism $\hA\to\ZZ_\ell$,
$v_{rj}\mapsto 0$.

\begin{lem}   \label{lem:triv}
 The following hold:
 \begin{itemize}
  \item[\rm(a)] Under the specialisation $\psi_1:\hA\to\ZZ_\ell$,
   $v_{rj}\mapsto 0$ (so $u_{rj}\mapsto\ze_{o(r)}^j$), the algebra $\cY$
   specialises to the group algebra of $TW$.
  \item[\rm(b)] The quotient of $\cY$ by the ideal $I$ generated by the
   $\{y_t-1\mid t\in T\}$ is isomorphic to the extension
   $\hA\otimes_A \cH(W,\bu)$ of the generic Hecke algebra of $W$.
 \item[\rm(c)] The natural $\hA$-module homomorphism $\cY\rightarrow \cH(W,\bu)$,
  $y_\br\mapsto h_\br$, in {\rm(b)} has a splitting
  $\cH(W,\bu)\to \hA[\ell^{-1}]\otimes_\hA \cY$ given by
  $h_\br\mapsto |T|^{-1}\sum_{t\in T}y_ty_\br$.
 \end{itemize}
\end{lem}

\begin{proof}
The first parts follows directly from the deformed order relation $(\dagger)$
and the corresponding result of Bessis \cite{Be01} for $B$. For (b), let $J_1$
be the ideal of $\hA[\hB]$ generated by the elements $\{t-1\mid t\in T\}$.
Then $I$ is the ideal of $\hA[\hB]$ generated by $J_1$ and the elements
$(\dagger)$ as $\br$ runs over
the braid reflections. Let $L$ be the ideal of $\hA[\hB]$ generated
by $J_1$ and the $\prod_{j=1}^{o(r)}(\br-u_{rj})$ as $\br$ runs over braid reflections. Then $\cH (W, \bu) = \hA[B] / L $. For an element 
$x= \prod_{j=1}^{o(r)}(\br-\zeta_{o(r)}^j(1 +E_r v_{rj}))\in J_1$ of the form
$(\dagger)$ set $x'=\prod_{j=1}^{o(r)}(\br-u_{rj} )$. Then 
$x +J_1=x' +J_1 $, whence $I =L $. 

For (c), note that $|T|^{-1}\sum_{t\in T}y_t$ is a central idempotent of
$\hA[\ell^{-1}]\otimes_\hA \cY$.
\end{proof}

We make some further straightforward observations. First, since for $\ell>2$
all reflections $r$ in an $\ell$-adic reflection group have order prime to
$\ell$, in that case $T_r\cong \ZZ/\ell^a\ZZ$.

For all braid reflections $\br$, the element $E_r$ commutes with $\br$ and
 $E_r^2=|T_r|E_r$.
Thus, over $\hA[\ell^{-1}]$, the element $E_r'=|T_r|^{-1}E_r$ is idempotent.
Multiplying ($\dagger$) with $E_r'$ respectively $1-E_r'$ we obtain the elements
$$(1-E_r')(\br^{o(r)}-1)\quad\text{and}\quad
  E_r'\prod_{j=1}^{o(r)}(\br-u_{rj}),\eqno(\dagger')$$
which generate the same ideal as ($\dagger$) over $\hA[\ell^{-1}]$. Thus,
($\dagger$)
``interpolates'' between the group relation for $r$ and the deformed Hecke
algebra relation ($\cH$) for $\br$. Let us also note the following:

\begin{lem}   \label{lem:T inv}
 The constant coefficient in the deformed order relation $(\dagger)$ is
 invertible in~$\hA[T]$.
\end{lem}

\begin{proof}
The constant term in the polynomial relation $(\dagger)$ for a braid reflection
$\br$ is (up to a root of unity) a product of factors $1+v_{rj}E_r$, which
has inverse $1-\zeta_{o(r)}^jv_{rj}/u_{rj}E_r\in \hA[T]$.
\end{proof}

We now give a more tangible description of $\cY(W,a,\bv)$. Recall that the
braid group $B$ has a presentation in terms of certain sets of braid
reflections together with so-called braid relations, encoded in \emph{braid
diagrams}, such that adding the order relations for the chosen braid
reflections, we obtain a presentation of $W$ (see \cite[Thm~2.27]{BMR} and
Bessis \cite[Thm~0.1]{Be01}). Choose reflections $r_1,\ldots,r_m$ in $W$
corresponding to a braid diagram for~$B$. It is known that any distinguished
reflection of $W$ is then conjugate to one of the $r_i$, and by their
construction all braid reflections projecting onto a fixed reflection of $W$
are conjugate in $B$. Then using Lemma~\ref{lem:T inv} we see that $\cY$ is the
associative unital $\hA$-algebra generated by elements
$\{y_t,y_{\br_i}\mid t\in T,\ 1\le i\le m\}$ subject to
\begin{itemize}
\item the $y_t$ satisfy the same relations as the corresponding group elements
 $t$ (i.e., they generate a subalgebra isomorphic to (possibly a quotient of)
 the group algebra $\hA[T]$);
\item the action relations between the $t,r_i$, with $t$ replaced by
  $y_t$ and $r_i$ by $y_{\br_i}$;
\item the braid relations between the $y_{\br_i}$; and
\item the deformed order relations $(\dagger)$ for the $y_{\br_i}$.
\end{itemize}

\subsection{On the structure of specialised Yokonuma type algebras.}
We show that under some additional hypothesis certain specialisations of $\cY$
are isomorphic to the group algebra of $TW$. The main results are
Theorem~\ref{thm:constr rep} and Theorem~\ref{thm:KOW2}.


For a specialisation $\psi:\hA \to R$ to a commutative ring $R$, let
$\cY_\psi:= R\otimes_\hA \cY$ denote the extension of scalars by $\psi$. Then
$\cY_\psi$ is the quotient of the group algebra $R\hB$ by the ideal
$$\Big\langle\prod_{j=1}^{o(r)}\big(\br-\ze_{o(r)}^j(1+ \psi(v_{rj})E_r)\big)
 \mid\text{$\br\in B$ braid reflection}\Big\rangle.$$

Let $W_0$ be a parabolic subgroup of $W$ and $B_0=B(W_0)$ be its braid group.
In \cite[\S2D]{BMR} is constructed an embedding $B_0\hookrightarrow B$,
well-defined up to $P$-conjugation, where $B=B(W)$, $P=P(W)$. By
\cite[Prop.~2.29 and~2.18]{BMR} this satisfies:

\begin{lem}   \label{lem:BMR}
 Let $W_0$ be a parabolic subgroup of $W$. Let $\tB_0$ be the inverse
 image of $W_0$ in $B$ and let $\tilde P_0$ be the subgroup of $P$ generated by
 the elements $\br^{o(r)}$, as $r$ runs over the distinguished reflections in
 $W\setminus W_0$. Let $B_0$ be the braid group of $W_0$. The above inclusion
 $B_0\hookrightarrow B$ has image contained in $\tB_0$ and induces an
 isomorphism $B_0\buildrel\sim\over\to \tB_0/\tilde P_0$.
\end{lem}

So we have the following diagram with exact columns:
$$\begin{matrix}
  1& & 1& & 1\\
  \downarrow& & \downarrow& & \downarrow\\
  P_0& \harr& P=\tilde P_0\rtimes P_0& =& P\\
  \downarrow& & \downarrow& & \downarrow\\
  B_0& \harr& \tB_0=\tilde P_0\rtimes B_0& \harr& B\\
  \downarrow& & \downarrow& & \downarrow\\
  W_0& =& W_0& \harr& W\\
  \downarrow& & \downarrow& & \downarrow\\
  1& & 1& & 1\\
\end{matrix}$$

\begin{rem}
Examples show that the isomorphism $B_0\buildrel\sim\over\to \tB_0/\tilde P_0$
might more generally hold for all reflection subgroups $W_0$ of $W$ generated
by distinguished reflections (see e.g. \cite[Prop.~3.24]{BMR}), so that the
assumptions of the subsequent Theorem~\ref{thm:constr rep} might be relaxed
accordingly. We will not need this here.
\end{rem}

If $W_\th$ is a parabolic subgroup of $W$, we will denote by $\bu_\th$ the set
$\bu_0$ in the notation of Section~\ref{subsec:braid} if $W_\th$ is the
reflection subgroup~$W_0$.

\begin{lem}   \label{lem:idempotent cutting}
 Let $R$ be an integral domain containing containing the $|T|$th roots of
 unity with field of fractions $K$ of characteristic $0$ and let $\psi:\hA\to R$
 be a specialisation. Let $I$ be the ideal of $R\hB$ generated by 
 $$\Big\{\prod_{j=1}^{o(r)}\big(\br-\ze_{o(r)}^j(1+ \psi(v_{rj})E_r)\big)
   \mid \text{$\br\in B$ braid reflection}\Big\}. $$
 Let $\th \in \Irr_K(T)$ and $e_\th\in KT$ the corresponding central
 idempotent. Suppose that the stabiliser $W_\th$ of $\th$ is a parabolic
 subgroup of $W$. Then
 $e_\th R\hB e_\th/e_\th I e_\th\cong\cH_\psi(W_\th,\bu_\th)$ as $R$-algebras.
 Here we regard $R\hB$ as a subset of $K\hB$. 
\end{lem}

Note that the assumption that $R$ contains the $|T|$th roots of unity is needed
in order to ensure that any ordinary irreducible character of $T$ is $R$-valued.

\begin{proof}
Let $\tB_\th$ be the full inverse image of $W_\th$ in $B$. For a braid
reflection $\br\in B$ set 
$$i_\br := \prod_{j=1}^{o(r)}\big(\br-\ze_{o(r)}^j(1+ \psi(v_{rj})E_r)\big)
  = (1- E'_r) (\br^{o(r)}-1) + E'_r \prod_{j=1}^{o(r)}(\br- \psi (u_{rj})) $$
where $ E_r'= |T_r|^{-1} E_r$. Then
$$\{e_\th x i_\br y e_\th
  \mid x,y\in \hB,\ \text{$\br\in B$ braid reflection}\}$$
generates $e_\th I e_\th$ as $R$-module. The set of braid reflections is
invariant under conjugation by~$B$, and $v_{rj} =v_{r'j}$ whenever $r$ and $r'$
are conjugate. Thus, if $x=t g$, $y= h s$ with $t,s \in T$ and $g, h\in B$ and
$\br$ is a braid reflection, then
$$ e_\th x i_\br y e_\th = \th (t) \th (s) e_\th i_{\tw{g}\br} gh e_\th,
  \quad \text{where }\th (t), \th(s) \in R. $$
Thus, $e_\th I e_\th $ is the $R$-span of $\{e_\th i_\br x e_\th\mid
x\in B,\ \br\in B\text{ braid reflection}\}$.

Now $r\in W_\th$ if and only if $\th(t^{-1} t^r) =1$ for all $t\in T$. Thus,
$e_\th E_r '= e_\th$ if $r\in W_\th$ and zero otherwise, and so
$$e_\th i_\br=\begin{cases} e_\th (\br^{o(r)}-1)& \text{if $\br\notin W_\th$,}\\
   e_\th\prod_j(\br-\psi(u_{rj}))& \text{if $\br\in W_\th$.}\end{cases}$$
Further, since $\br^{o(r)}\in P$ commutes with $T$, we have 
$e_\th(\br^{o(r)}-1) x e_\th = (\br^{o(r) }-1)e_\th x e_\th$
and if $r\in W_\th$, then 
$$e_\th\prod_j(\br-\psi(u_{rj}))x e_\th=\prod_j(\br-\psi(u_{rj}))e_\th xe_\th.$$
For any $x\in B$, $e_\th x e_\th= xe_\th $ if $x\in\tB_\th$ and
zero otherwise. Hence $e_\th I e_\th$ is the $R$-span of
$$\big\{(\br^{o(r)} -1 )x e_\th\mid x\in\tB_\th,\ r\notin W_\th\big\}\cup
 \big\{\prod_j(\br -\psi (u_{rj}))x e_\th\mid x\in\tB_\th,\ r\in W_\th\big\}.$$
By the same argument we have that $e_\th R\hB e_\th = e_\th R [T\tB_\th] e_\th$
and since $e_\th $ is $T\tB_\th $-stable and $e_\th$ is idempotent we also have
$ e_\th R [ T\tB_\th ] e_\th = R [T\tB_\th ] e_\th $.
Since $\th$ is linear and $\tB_\th$-stable, there is an $R$-algebra isomorphism
$R\tB_\th \cong R[T\tB_\th]e_\th$ given by $x\mapsto x e_\th$. This induces an
isomorphism
$$R\tB_\th/J\cong R [T \tB_\th] e_\th /e_\th I e_\th$$
where $J\unlhd R\tB_\th$ is the ideal generated by $\{\br^{o(r)}-1\mid
r\notin W_\th\}\cup\{\prod_j (\br -\psi(u_{rj}))\mid r\in W_\th\}$. By
Lemmas~\ref{lem:BMR} and~\ref{lem:triv}(b),
$R\tB_\th/J\cong \cH_\psi(W_\th,\bu_\th)$. 
\end{proof}

\begin{thm} \label{thm:constr rep}
 Assume all stabilisers $W_\th$ of elements $\th\in\Irr(T)$ are parabolic
 subgroups of $W$. Let $K$ be a field of characteristic $0$ containing the
 $|T|$th roots of unity and let $\psi:\hA\to K$ be a ring homomorphism. Then
 \begin{enumerate}
  \item [(a)]
   $$\cY_\psi \cong\prod_\th \Mat_{|W: W_\th|}(K)\otimes_K \cH_\psi(W_\th,\bu_\th)$$
   as $\th$ runs over a set of representatives of $W$-orbits on $\Irr_K(T)$;
  \item[(b)] $\dim_K \cY_\psi =|TW|$.
  \item[(c)] Suppose that $K_W\subseteq K $ and $\psi$ is the inclusion
   homomorphisms. Then $\cY_\psi\cong K[TW]$.
   \end{enumerate}
\end{thm}
 
\begin{proof}
Part~(a) is immediate from Lemmas~\ref{lem:infclifford}
and~\ref{lem:idempotent cutting} applied with $R= K$. Part~(b) follows from
part~(a) by Theorem~\ref{thm:Hecke free} and Lemma~\ref{lem:infclifford}(b)
applied with $G = T W$.   \par
Now assume $K$ and $\psi$ are as in~(c). As explained in
Section~\ref{subsec:spec}, for all $\th\in\Irr(T)$,
$\cH_\psi(W_\th,\bu_\th) = K\otimes_\hA \cH(W_\th,\bu_\th)\cong KW_\th$.
Now~(c) follows from~(a) and Lemma~\ref{lem:infclifford}(a) applied with
$G=T W$ and $I=0$, noting that $K (TW)_\th e_\th \cong KW_\th$.
\end{proof}
 
We now turn to specialisations of $\hA$ to finite extensions of $\ZZ_\ell$. For
the rest of this subsection, the following notation will be in effect. Let
$\ZZ_\ell\subseteq\cO$ be a complete discrete valuation ring with uniformiser
$\pi$, residue field $k$ and field of fractions $K$. Let $\psi:\hA\to\cO$ be a
$\ZZ_\ell$-algebra homomorphism and denote by $\bar\psi$ the composition of
$\psi$ with the canonical map $\cO\to k$. Recall that for $x\in\hA[\hB]$ we
denote by $y_x$ its image in $\cY$. For $y\in\cY$ let
$\tilde y:=1_\cO\otimes y\in\cY_\psi$, and
$\bar y:=1_k\otimes y\in \cY_{\bar \psi}$.

\begin{lem}   \label{l:fg}
 Let $\bW$ be a set of coset representatives of $P$ in $B$ and let~$J$ be the
 ideal of $\cY_\psi$ generated by $\{\tilde y_t-1 \mid t\in T\}$ and $\pi$.
 Then, $kW\cong\cY_{\psi}/J$ via the map which sends $w\in W$ to
 $\tilde y_\bw +J$ for $\bw\in\bW$ lifting $w$.
\end{lem} 

\begin{proof}
For a distinguished reflection $r\in W$, $\ell$ dvides $|T_r|$, hence
$u_{rj}\in I$. Now the result follows from Lemma~\ref{lem:triv}(b) (suitably
adapted to the coefficient ring $\cO $).
\end{proof} 

\begin{thm}   \label{thm:O finite generation}
 Let $\bW$ be a set of coset representatives of $P$ in $B$. If all stabilisers
 $W_\th$ of elements $\th\in \Irr(T)$ are parabolic subgroups of $W$, then 
 $X:=\{\tilde y_t \tilde y_{\bw} \mid t\in T,\ \bw \in\bW\}$ is an $\cO$-basis
 of~$\cY_\psi$.
\end{thm} 

\begin{proof}
We first show that $\{\bar y_t\bar y_{\bw}\mid t\in T,\ \bw \in\bW\}$ generates
$\cY_{\bar \psi}$ as $k$-vector space. Let $R \subseteq \cY_{\bar \psi}$ be the
$k$-span of $\{\bar y_t\mid t\in T\}$.
Then $R$ is a commutative $k$-subalgebra of $\cY_{\bar \psi}$. Let $Q$ be the
ideal of $R$ generated by $\{\bar y_t -1\mid t\in T\}$. Since $T$ is a finite
abelian $\ell$-group, $Q$ is a nilpotent ideal of $R$. We consider
$\cY_{\bar \psi} $ as a left $R$-module. Since $T$ is normal in $\hB$, the
$R$-submodule $Q\cY_{\bar \psi}$ of $\cY_{\bar \psi}$ is an ideal of
$\cY_{\bar\psi}$. Further, $Q\cY_{\bar\psi}$ is the image of the ideal $J$ of
Lemma~\ref{l:fg} in $\cY_\psi/\pi\cY_\psi\cong \cY_{\bar\psi}$ and for any
$\bw \in \bW$, $\bar y_\bw + Q\cY_{\bar \psi} = \tilde y_{\bw} + J$. So, by
Lemma~\ref{l:fg}, $\{\bar y_{\bw} + Q\cY_{\bar\psi}\mid \bw\in\bW\}$
generates $\cY_{\bar \psi}/Q\cY_{\bar \psi}$ as $k$-vector space and hence as
$R$-module. Applying Lemma~\ref{lem:Nakayama} with $M = \cY_{\bar\psi}$ and the
nilpotent ideal $Q$ we obtain that $\{\bar y_\bw \mid \bw\in\bW\}$ generates
$\cY_{\bar\psi}$ as $R$-module. Since $R$ is generated by
$\{\bar y_t\mid t\in T\}$ as $k$-vector space we have the required result.

Now we claim that in order to prove the theorem it suffices to prove that
$\cY_\psi$ is finitely generated as $\cO$-module. Indeed, suppose that
$\cY_\psi$ is finitely generated as $\cO$-module. Then by the previous
paragraph  and the standard Nakayama lemma (Lemma~\ref{lem:Nakayama}) applied
to the ring $\cO$ and ideal $\pi\cO$, the set
$X=\{\tilde y_t \tilde y_{\bw} \mid t\in T, \bw \in\bW\}$ generates $\cY_\psi$
as $\cO$-module. Then $1\otimes X$ generates $K \otimes_\cO\cY_\psi$ as
$K$-vector space and since the latter has dimension $|TW|$ by
Theorem~\ref{thm:constr rep}(b), $X$ is an $\cO$ basis of $\cY_{\psi}$. 

It remains only to show that $\cY_\psi$ is finitely generated as $\cO$-module.
For this, suppose first that $\cO$-contains the $|T|$th roots of unity. Let $I$
be the ideal of $\cO\hB$ generated by the $(\dagger)$-relations, let 
$$I' = \bigoplus_{\th, \mu \in \Irr(T)} e_\th I e_\mu
  \subseteq \bigoplus_{\th, \mu\in \Irr(T)} e_\th \cO \hB e_\mu $$
with $e_\th\in KT$ as in Lemma~{\rm\ref{lem:idempotent cutting}} and let
$\tilde I = \cO \hB \cap I'$. Since $I'$ is an ideal of
$\bigoplus e_\th\cO\hB e_\mu$, $\tilde I$ is an ideal of
$\cO\hB$ containing $I$. On the other hand, $|T|e_\th\in\cO\hB$ for all $\th$,
hence $|T|^2 e_\th Ie_\mu\subseteq I$ and $|T|^2\tilde I\subseteq I$.
The kernel of the composition of the inclusion $\cO\hB\hookrightarrow
\bigoplus e_\th\cO \hB e_\mu$ with the surjection
$\bigoplus e_\th \cO \hB e_\mu\twoheadrightarrow\bigoplus e_\th\cO\hB e_\mu/ I'$
is $\tilde I$. Thus, $\cO \hB /\tilde I $ is isomorphic to a submodule of
$\bigoplus e_\th \cO \hB e_\mu / I'$. On the other hand,
$$\bigoplus_{\th,\mu\in\Irr(T)} e_\th \cO \hB e_\mu / I' \cong
  \bigoplus_{\th,\mu \in\Irr(T)} e_\th \cO \hB e_\mu / e_\th I e_\mu. $$
If $\mu = \,^x \th $ for $x \in B$, then
$ e_\th \cO \hB e_\mu/ e_\th I e_\mu \cong e_\th \cO \hB e_\th /e_\th I e_\th$
via right multiplication by $x$ and it follows from
Lemma~\ref{lem:idempotent cutting} and Theorem~\ref{thm:Hecke free} that
$e_\th\cO \hB e_\th / e_\th I e_\th$ is finitely generated free as $\cO$-module.
If $\th, \mu$ are in different $W$-orbits, then $ e_\th \cO \hB e_\mu = 0$. By
the above displayed equation, $\bigoplus e_\th \cO \hB e_\mu / I'$ is finitely
generated free as $\cO$-module. Since $\cO $ is a principal ideal domain, and
since $\cO \hB /\tilde I $ is isomorphic to a submodule of
$\bigoplus e_\th \cO \hB e_\mu / I'$, it follows that $\cO \hB /\tilde I $ is
finitely generated free as $\cO$-module. Also, $\pi^r(\tilde I/I ) = 0$
for $r$ equal to twice the $\pi$-adic valuation of $|T|$. We saw above that
$\cY_\psi/\pi \cY_\psi \cong \cY_{\bar\psi}$ is finitely generated as $k$-vector
space and hence as $\cO$-module. Thus by Lemma~\ref{lem:lifting} applied
with $M= \cY_\psi= \cO \hB /I $ and $N=\tilde I/I $ we have $\cY_\psi$ is
finitely generated as $\cO$-module.

Now consider the general case. Let $\cO'$ be a finite extension of $\cO$
containing the $|T|$th roots of unity and let $\psi'$ be the composition of
$\psi$ with inclusion of $\cO$ in $\cO'$. By the previous part, applied with
$\psi'$ in place of $\psi$, $\cY_{\psi'}$ is finitely generated as
$\cO'$-module. Since $\cO'$ is a finite extension of $\cO$, $\cY_{\psi'}$ is
also finitely generated as $\cO$-module. Since $\cO$ is a direct summand of $\cO'$ as 
$\cO$-module, the inclusion
$\cO \hookrightarrow \cO'$ is pure. Thus the map
$\cY_\psi\to\cO'\otimes_\cO\cY_\psi\cong \cY_{\psi'}$, $y\mapsto 1\otimes y$,
is injective and consequently $\cY_\psi$ is isomorphic to an $\cO$-submodule of
$\cY_{\psi'}$. Since $\cO$ is Noetherian and since as shown above $\cY_{\psi'}$
is finitely generated as $\cO$-module, $\cY_\psi$ is finitely generated as
$\cO$-module. 
\end{proof}

The following is an application of a theorem of K\"ulshammer, Okuyama and
Watanabe (see \cite[Thm~4.8.2]{Li18}).
Recall that if $R$ is a commutative ring and $C$ is a subalgebra of an
$R$-algebra $B$, then $B$ is \emph{relatively $C$-separable} if $B$ is a direct
summand of $B\otimes_C B$ as a $(B,B)$-bimodule. 

\begin{thm}   \label{thm:KOW2}
 Suppose that $W$ is an $\ell'$-group. Then there exists an $\cO$-algebra
 isomorphism $\cY_\psi\cong \cO[TW]$ sending $\tilde y_t$ to $t$ for any
 $t\in T$.
\end{thm} 

\begin{proof}
Let $\bW$ be a set of coset representatives of $P$ in $B$.
Since $W$ is an $\ell'$-group, $W_\th$ is a parabolic subgroup of $W$ for all
$\th\in\Irr(T)$. Thus, by Theorem~\ref{thm:O finite generation},
$\{\tilde y_t \tilde y_\bw\mid t\in T,\ \bw\in\bW\}$ is an $\cO$-basis of
$\cY_\psi$. Further, we may regard $\cO T$ as an $\cO$-subalgebra of $\cY_\psi$
via the identification of $\cO T$ with the subalgebra generated by
$\{\tilde y_t\mid t\in T\}$. Under this identification, again via
Theorem~\ref{thm:O finite generation}, there is a homomorphism of
$(\cO T,\cO T)$-bimodules $\gamma:\cO[TW]\to \cY_\psi$ defined by
$\gamma(tw)=\tilde y_t \tilde y_\bw$. 
 
Let $J$ be the ideal of $\cY_\psi$ generated by $\{\tilde y_t-1\mid t\in T\}$
and $\pi$. It follows from Lemma~\ref{l:fg} that the composition of $\gamma$
with the natural surjection $\cY_\psi\to \cY_\psi/J$ is an $\cO$-algebra
homomorphism. Now we may apply \cite[Thm~4.8.2]{Li18} to obtain an
$\cO$-algebra homomorphism $\sigma: \cO[TW] \to \cY_\psi$ extending the
$\cO$-algebra homomorphism $\cO[TW]\to \cY_\psi/J$ obtained above from $\gamma$
and satisfying $\sigma (t) =\tilde y_t$ for all $t\in T$. For this one needs to
have that $J$ is contained in the radical of $\cY_\psi$ and that $\cO [TW]$ is
relatively $\cO T$-separable. The second condition holds since $W$ is an
$\ell'$-group (see \cite[Prop.~2.6.9]{Li18}) whereas the first condition holds
since $T$ is a finite normal $\ell$-subgroup of $\hB$ and, by
Theorem~\ref{thm:O finite generation}, $\cY_\psi$ is finitely generated as
$\cO$-module. 

The surjectivity of $\sigma $ follows by Nakayama's lemma since the composition
of $\sigma$ with $\cY_\psi\to\cY_\psi/I$ is surjective and then the injectivity
follows since both algebras are free of the same rank.
\end{proof}

\begin{rem}   \label{rem:stab}
The assumption of Theorem~\ref{thm:constr rep} on stabilisers is satisfied
whenever $\ell$ is very good for $W$, e.g., when $|W|$ is coprime to $\ell$,
or if $W=G(e,1,n)$ with $e\ge2$, see \cite[Prop.~2.3]{KMS}.
\end{rem}

\subsection{Freeness} \label{sub:free}
We propose the following, analogous to the (now proven) Freeness Conjecture
(Theorem~\ref{thm:Hecke free}) for cyclotomic Hecke algebras:

\begin{conj}   \label{conj:free}
 The algebra $\cY$ is free over $\hA$ of rank $|TW|$. More precisely, 
 there is a section $W\to\bW\subset B$ of the natural map $B\to W$
 containing~$1$ such that $\{y_ty_\bw\mid t\in T,\,\bw\in\bW\}$ is an
 $\hA$-basis of $\cY$.
\end{conj}

For Weyl groups and parameters occurring in finite reductive groups, the
freeness follows from the construction as an endomorphism algebra, and the
dimension from the number of double cosets of a maximal unipotent subgroup,
that is, the Bruhat decomposition; see Lusztig \cite[34.2--34.10]{Lu7} for a
detailed investigation.
We propose a proof in the case of finite Coxeter groups and for most infinite
series of complex reflection groups.

\begin{thm}   \label{thm:free Cox}
 Conjecture~{\rm\ref{conj:free}} holds for any finite Coxeter group.
\end{thm}

\begin{proof}
Assume that $W$ is a Coxeter group and choose a presentation of $B$ on braid
reflections $\br_1,\ldots,\br_m\in B$ mapping to the Coxeter generators
of~$W$. Clearly, the set of all monomials in the $y_t,y_{\br_i}$ forms
a generating system for $\cY$ as an $\hA$-module. By the `action relations' any
such monomial can be rewritten into an $\hA$-linear combination of elements
$y_ty_\bw$ with $t\in T$ and $\bw$ a monomial in the generators $\br_i$,
$1\le i\le m$. Now by Matsumoto's lemma, by using the braid relations plus
the quadratic relations $(\dagger)$ expressing $y_{\br_i}^2$ as a linear
combination of smaller powers of $y_{\br_i}$, $\bw$ can be rewritten into an
$\hA[T]$-linear combination of elements from a fixed set $\bW\subset B$ of
reduced expressions of elements of~$W$.   \par
Thus any monomial in the generators is an $\hA$-linear combination of elements
$y_ty_\bw$ with $t\in T$ and $\bw\in\bW$. Since the $y_t$ satisfy the
same relations as the corresponding $t\in T$, there are at most $|T|$ distinct
elements $y_t$, so we have identified a generating system for
$\cY$ of cardinality $|TW|$. By Theorem~\ref{thm:constr rep} this must be free
over $K$, hence an $\hA$-basis of~$\cY$.
\end{proof}

\begin{thm}   \label{thm:free G(e,p,n)}
 Conjecture~{\rm\ref{conj:free}} holds for $W=G(e,p,n)$ with $e|(\ell-1)$ for
 any divisor $p$ of $e$, except possibly when $n=2$, $e,p$ are both even and
 $p\ne e$.
\end{thm}

\begin{proof}
The group $W=G(e,p,n)$ is a normal reflection subgroup of $W_1:=G(e,1,n)$ of
index~$p$. First assume that $p<e$. Then the braid group $B$ of $W$ is normal
in the braid group $B_1$ of~$W_1$ of index~$p$ by \cite[\S3.B1]{BMR}. Also,
the corresponding tori $T$ can be identified, such that $T.B$ is normal in
$T.B_1$ of index~$p$. A system of coset representatives is given by
$\{\br_1^i\mid 0\le i\le p-1\}$, where $\br_1\in B_1$ lifts a distinguished
reflection $r\in W_1$ of order~$e$. Let $\zeta_e$ be a primitive $e$th root of
unity and set $p':=e/p$. Recall the parameters $u_{rj}$, $1\le j\le p'$, for
$\cY$ at the reflection $r$. Let $K$ be a sufficiently large extension of
$\Frac(\hA)$. Consider parameters
$u_{rj}':=u_{rj}^{1/p}$ for $1\le j\le p'$, and
$u_{r,j+p'}':=\zeta_e^{p'} u_{rj}'$ for $1\le j\le e-p'$. Now over $K$, the
relation~($\dagger$) for $y_1:=y_{\br_1}$ can be rewritten as
$$(1-E')(y_1^e-1)+E'\prod_{j=1}^e(y_1-u_{rj}')=0,$$
with $E':=|T_r|^{-1}E_r$ (see ($\dagger'$) above). Note that
$\prod_{i=0}^{p-1}(y_1-u_{r,j+p'i}')=y_1^p-u_{rj}$ for any $j$. Similarly,
over $K$ the relation~($\dagger$) for the generator $y_1^p$ of $\cY$ can be
written as
$$(1-E')(y_1^{pp'}-1)+E'\prod_{j=1}^{p'}(y_1^p-u_{rj})=0.$$
Thus, we obtain the same relation for $y_1^p$ in $\cY$ as before.
\par
Let $\cY_1$ be the quotient of $\hA[T.B_1]$ by the deformed order relations,
which agree with those for $\cY$ as we just saw and hence can be written over
$\hA$, except in the excluded case $n=2$, $e,p$ both even, when $G(e,p,2)$
contains an additional class of reflections.   \par
We claim that the conjecture holds for $W_1$. Our proof for this
closely follows some arguments in Bremke--Malle \cite{BM97}. Let
$\br_1,\ldots,\br_n\in B$ be braid reflections corresponding to the standard
presentation, so that $\br_2,\ldots,\br_n$ generate the braid group on $n$
strands (of type $A_{n-1}$) and $(\br_1\br_2)^2=(\br_2\br_1)^2$. Set
$y_2:=y_{\br_2}$. Now by Lemma~\ref{lem:rel G(e,1,n)} below for any $a,b\ge1$ we
have
$$y_2y_1^ay_2y_1^b=\al y_1^by_2y_1^ay_2+
  \sum_{i=1}^b(\al_iy_1^{a+b-i}y_2y_1^i+\al_i'y_1^iy_2y_1^{a+b-i})$$
for suitable $\al\in \hA[T]^\times$ and $\al_i,\al_i'\in \hA[T]$. With this,
one deduces as in the proof of \cite[Prop.~2.4]{BM97} that there is a set
$\cB_1\subset B_1$ of cardinality $|W_1|$ consisting of monomials in the
$\br_i$, as in \cite[Lemma~1.5]{BM97}, such that any monomial in the $y_t,y_i$
can be
rewritten in $\cY_1$ into an $\hA$-linear combination of the $|TW_1|$ products
$\cB:=\{y_ty_\bw\mid t\in T,\ \bw\in\cB_1\}$. Thus, $\cB$ is linearly
independent over $K$ by Theorem~\ref{thm:constr rep} and so an $\hA$-basis
of $\cY_1$. This proves our claim for $W_1$.
\par
Now $\hA[T.B_1]=\bigoplus_{i=0}^{p-1}\hA[T.B]\br_1^i$ is $\ZZ/p\ZZ$-graded and
multiplication by $\br_1$ defines $\hA$-module isomorphisms between the
summands. Furthermore, the defining ideal $I$ for $\cY$ in $\hA[T.B]$ is
contained in the defining ideal $I_1$ of $\cY_1$ in $\hA[T.B_1]$, and
$I_1=\bigoplus_{i=0}^{p-1}I\br_1^i$ is graded. So
$\cY_1=\hA[T.B_1]/I_1=\bigoplus_{i=0}^{p-1}\cY y_1^i$ and multiplication with
$y_1$ induces $\hA$-module isomorphisms between the summands on the right.
By construction the $\hA$-basis $\cB$ of $\cY_1$ has the property that
$\cB=\bigcup_{i=0}^{p-1}(\cB\cap \cY y_1^i)$, hence $\cB\cap \cY$ is an
$\hA$-free generating system of~$\cY$.   \par
Finally assume that $p=e$. Since $G(e,e,2)$ is a Coxeter group (the dihedral
group of order $2e$), by Theorem~\ref{thm:free Cox} we may assume $n\ge3$. In
this case, the braid group $B$ of $W=G(e,e,n)$ is a normal subgroup of index
$e$ of the quotient $\bar B_1=B_1/\langle \br_1^e\rangle$ of the braid group of
$W_1$ (see \cite[Prop.~3.24]{BMR}). Thus, $\cY$ is an $\hA$-subalgebra of
$\hA[T.B_1]/I$ where $I$ is generated by $\br_1^e-1$ and the relations
$(\dagger$) for $\br_2,\ldots,\br_n$. We can now argue precisely as in the
previous case.
\end{proof}

The following was used in the preceding proof:

\begin{lem}   \label{lem:rel G(e,1,n)}
 Let $\cY=\cY(G(e,1,n))$ and $y_1,y_2\in \cY$ images of braid reflections
 satisfying $y_2y_1y_2y_1=y_1y_2y_1y_2$ and such that the
 corresponding reflections $r_1,r_2\in W$ have order~$e,2$ respectively. Then
 for all integers $a,b\ge1$ there exist $\al\in \hA[T]^\times$ and
 $\al_i,\al_i'\in \hA[T]$ such that
 $$y_2y_1^ay_2y_1^b
  =\al y_1^by_2y_1^ay_2
  +\sum_{i=1}^b(\al_iy_1^{a+b-i}y_2y_1^i+\al_i'y_1^iy_2y_1^{a+b-i}).$$
\end{lem}

\begin{proof}
Write the relation~($\dagger$) for $y_2$ as $y_2^2=\la y_2+\mu$ with
$\mu\in \hA[T]^\times$ and $\la\in \hA[T]$, so
$y_2^{-1}=\mu^{-1}y_2-\la\mu^{-1}$. The relation between $y_1,y_2$
implies $y_1^ay_2y_1y_2=y_2y_1y_2y_1^a$ for all $a\ge1$. Thus
we find
$$\begin{aligned}
  y_2y_1^ay_2y_1=&\ y_2y_1^ay_2y_1\cdot y_2y_2^{-1}=y_2^2y_1y_2y_1^ay_2^{-1}
    =(\la y_2+\mu)y_1y_2y_1^a(\mu^{-1}y_2-\la\mu^{-1})\\
  =&\ \mu y_1y_2y_1^ay_2\mu^{-1}+\la y_2y_1y_2y_1^ay_2^{-1}-\mu y_1y_2y_1^a\la\mu^{-1}\\
  =&\ \mu y_1y_2y_1^ay_2\mu^{-1}+\la y_1^ay_2y_1-\mu y_1y_2y_1^a\la\mu^{-1}\\
  =&\ \mu'y_1y_2y_1^ay_2+\la y_1^ay_2y_1-\la'y_1y_2y_1^a
\end{aligned}$$
for suitable $\mu'\in \hA[T]^\times$ and $\la'\in \hA[T]$, giving the claim for
$b=1$. For $b=2$, using the previous result twice we find
$$\begin{aligned}
  y_2y_1^ay_2y_1^2
  =&\ (\mu'y_1y_2y_1^ay_2+\la y_1^ay_2y_1-\la'y_1y_2y_1^a)y_1\\
  =&\ \mu'y_1y_2y_1^ay_2y_1+\la y_1^ay_2y_1^2-\la' y_1y_2y_1^{a+1}\\
  =&\ \mu'y_1(\mu'y_1y_2y_1^ay_2+\la y_1^ay_2y_1-\la'y_1y_2y_1^a)+\la y_1^ay_2y_1^2-\la' y_1y_2y_1^{a+1}\\
  =&\ \mu''y_1^2y_2y_1^ay_2+\sum_{i=1}^2(\al_iy_1^{a+2-i}y_2y_1^i+\al_i'y_1^iy_2y_1^{a+2-i})
\end{aligned}$$
for suitable $\mu'',\al_i,\al_i'$. A straightforward induction yields the claim
for arbitrary $b$.
\end{proof}

\begin{rem}   \label{rem:Marin}
Marin \cite[Def.~5.4]{Ma18a} defines for arbitrary complex reflection groups
$W$ an $\hA$-algebra $M$ attached to $W$ as follows: let $\cL$ be the lattice
of intersections of the hyperplane arrangement of $W$. Then $M$ is the quotient
of the group algebra over $\hat A$ of the semidirect product $\cL\rtimes B(W)$
by the deformed order relations $(\dagger)$ for the braid reflections of
$B(W)$. This algebra is generated by images of braid reflections $\br'$ and
idempotents $e_r$, $r\in W$ a reflection (by \cite[5.1]{Ma18a}).
Marin shows \cite[Thm~1.3]{Ma18b} that $M$ is a free $\hA$-module of
rank $|W||\cL|$. If $W$ is $\ell$-adic, there is a natural morphism
$$i_W:M\to \hA[\ell^{-1}]\otimes_\hA \cY,\quad
  \br'\mapsto y_\br,\ e_r\mapsto\ell^{-a}E_r,$$
from Marin's algebra to ours. We expect this to be injective, but in general
far from surjective, since his algebra is free of rank independent of $\ell$
(compare to Theorem~\ref{thm:constr rep}).
\end{rem}

\subsection{A trace form}   \label{subsec:trace form}
Assume for the rest of the section that $W_\th$ is a parabolic subgroup of
$W$ for all $\th\in\Irr(T)$; this holds whenever $\ell$ is very good for
$(W,L)$, see Remark~\ref{rem:stab}. Let $K$ be an extension of $\QQ_\ell$ by
the $\ell^a$th roots of unity. Let $\tilde\bu=(\tu_{rj})$ be as in
Section~\ref{subsec:braid} and let $\tK=\Frac(K[\tilde\bu])$. Recall from
Section~\ref{subsec:spec} that for any $\th\in\Irr(T)$,
$\tK\otimes\cH(W_\th,\bu_\th)\cong \tK W_\th$ is split semisimple and the
irreducible characters of $\cH(W_\th,\bu_\th)$ over $\tK$ are identified with
$\Irr(W_\th)$. Here, as before we denote by $\bu_\th$ the set $\bu_0$ in the
notation of
Section~\ref{subsec:braid} if $W_\th$ is the reflection subgroup~$W_0$. Then,
with $U$ a $\tK\otimes\cH(W_\th,\bu_\th)$-module affording $\phi$ and 
$U_{\th,\phi}:=\Ind_{\hB_\th}^\hB(U)$, Theorem~\ref{thm:constr rep}(a) shows
$$\Irr(\cY_\tK )=\{U_{\th,\phi}\mid \th\in \Irr(T)/W,\ \phi\in\Irr(W_\th)\}.$$
We let $\chi_{\th,\phi}$ denote the character of $U_{\th,\phi}$. 

We consider the following non-degenerate trace form $\cY \to \tK$:
\begin{equation}   \label{eqn:form}
 \tau:=\tau_\cY
 :=\sum_{\th/W}\sum_{\phi\in\Irr(W_\th)}\frac{1}{f_{\th,\phi}}\chi_{\th,\phi}.
\end{equation} 
Here, for any $\th\in\Irr(T)$, $f_{\th,\phi} \in \tilde A $ is the Schur
element of the Hecke algebra $\cH(W_\th,\bu_\th)$ indexed by $\phi$ as in
Section~\ref{subsec:braid}.

\begin{prop}   \label{prop:is form}
 Assume that $W_\th$ is a parabolic subgroup of $W$ for all $\th\in\Irr(T)$.
 Assume also that $\cH(W,\bu)$ is strongly symmetric with respect to
 $\bW\subset B$ as in Definition~{\rm\ref{def:strongly quasi}}. Then we have
 $$\tau(y_ty_\bw)=\delta_{t,1}\delta_{\bw,1}|T|
  \qquad\text{for any $t\in T$, $\bw\in\bW$}.$$
\end{prop}

\begin{proof}
For $\th\in\Irr(T)$ let $C_\th$ be a system of coset representatives of
$W_\th$ in $W$, and $\bC_\th\subseteq\bW$ the corresponding system of coset
representatives of $T.\tB_\th$ in $T.B$. Let $\th\in\Irr(T)$ and
$\phi\in\Irr(W_\th)$, and let $U$ be a corresponding representation of
$\cH(W_\th,\bu_\th)$ (which we consider as a representation of
$\hB_\th =T\tB_\th$ as above). Let us set $U^0(x):=U(x)$ if
$x\in T.\tB_\th$ and $0$ otherwise. Then
$$U_{\th,\phi}(y_ty_\bw)=\sum_{\bx\in\bC_\th} U^0((y_ty_\bw)^\bx)
  =\sum_{\bx\in\bC_\th,y_\bw^\bx\in\tB_\th}\th(y_t^\bx)U(y_\bw^\bx),$$
so
$$\sum_{\phi\in\Irr(W_\th)}\frac{1}{f_{\th,\phi}}\chi_{\th,\phi}(y_ty_\bw)
  =\sum_{\bx\in\bC_\th,y_\bw^\bx\in\tB_\th}\th(y_t^\bx)\sum_{\phi\in\Irr(W_\th)}\frac{1}{f_{\th,\phi}}\chi_{\th,\phi}(y_\bw^\bx)
  =\!\!\!\!\sum_{\bx\in\bC_\th,y_\bw^\bx\in\tB_\th}\th(y_t^\bx)t_{W_\th,\bu_\th}(y_\bw^\bx)$$
with $t_{W_\th,\bu_\th}$ as in Definition~\ref{def:strongly quasi}. By the
choice of $\bW$ we have
$$t_{W_\th,\bu_\th}(y_\bw^\bx)=t_{W,\bu}(y_\bw^\bx)=t_{W,\bu}(y_\bw)
  =\delta_{\bw,1}.$$
Thus the form $\tau$ evaluates to
$$\tau(y_ty_\bw)
  =\sum_{\th/\cF}\sum_{\bx\in\bC_\th}\th(y_t^\bx)\delta_{\bw,1}
  =\sum_{\th\in\Irr(\bT)}\th(y_t)\delta_{\bw,1}
  =\delta_{t,1}\delta_{\bw,1}|T|,$$
as desired.
\end{proof}

It seems natural to ask the following:

\begin{ques}   \label{conj:trace}
 Let $(W,L)$ be a simply connected $\ell$-adic reflection group for which
 $\ell$ is very good. Does the form $|T|^{-1}\tau$ take values in $\hat A$ and
 is it then a symmetrising form on~$\cY$ over~$\hat A$?
\end{ques}


Note that an affirmative answer to the first part of Question~\ref{conj:trace}
follows under the assumptions of Proposition~\ref{prop:is form}.

\subsection{Relation to classical Yokonuma algebras}   \label{sub:classical yokonuma}
Suppose that $W$ is the Weyl group with respect to a maximally split torus
$\bT_0$ of a connected reductive group $\bG$ with an $\FF_q$-structure defined
by a split Frobenius map $F:\bG\to\bG$. Set $G =\bG^F$ and set $T_0 =\bT_0^F$.
Let $\bU$ be the unipotent radical of an $F$-stable Borel subgroup of $\bG$
containing $\bT_0$, and let $U =\bU^F$. Let $\ell$ be a prime dividing $q-1$ 
and set $e_U =|U|^{-1}\sum_{u\in U}u \in \ZZ_\ell G$. Then 
$$\cY':=\End_{\ZZ_\ell G}(\ZZ_\ell [G/U ]) = e_U\ZZ_\ell G e_U $$
is the associated classical Yokonuma Hecke algebra \cite{Yo67}.

Let $r_1,\ldots,r_m$ be Coxeter generators of the Weyl group $W$. For
$t\in T_0$, let $t'=e_U te_U$ and for each $i$,
let $E_i:= E_{[T_0, r_i] }e_U$, where for any subgroup $A\le G$ we denote by
$E_A$ the sum of elements of $A$. By \cite[Thm~2]{JK01}, $\cY'$ has a
generating set $\{t',s_i\mid t\in T_0,\ 1\le i\le m\}$ such that 
\begin{itemize}
\item the $t'$ satisfy the same relations as the corresponding group elements
 $t$;
\item the action relations between the $t, r_i$, with $t$ replaced by
 $t'$ and $r_i$ by $s_i$ hold in $\cY'$;
\item the braid relations between the $r_i$, with $r_i$ replaced by $s_i$ hold
 in $\cY'$; and
\item $s_i^2 = 1- q^{-1} (E_i - s_i E_i)$.
\end{itemize}
Note that Juyumaya--Kannan work over the complex numbers but it can be checked
from the explicit description of the $s_i$s in terms of the standard
generators coming from the Bruhat decomposition, that the above holds over any
ring in which $q$ is invertible.

\begin{prop}   \label{prop:classical yokonuma}
 Suppose that there is a $W$-equivariant isomorphism between $T$ and the Sylow
 $\ell$-subgroup of $T_0$. Let $H$ be the $\ell'$-Hall subgroup of $T_0$, let
 $e_{H} =|H|^{-1}E_H$ be the principal block idempotent of $\ZZ_\ell T_0$ and
 set $f=e_{H} e_U$. Let $\psi_q:\hA \to\ZZ_\ell$ be the specialisation
 corresponding to $u_{\br1}\mapsto -1$, $u_{\br2}\mapsto q$ for all~$\br$. Then
 there is an isomorphism of $\ZZ_\ell$-algebras $\cY_{\psi_q}\cong f \cY' f$.
\end{prop} 

\begin{proof}
Note that $f$ is an idempotent of $\ZZ_\ell G$. By considering the generating
set of $\cY'$ described above, one sees that $f$ is central in $\cY'$ and for
any $x\in H$, $xf = f$, hence
$\{t'f, s_if\mid t\in T_0,\ 1\le i\le m\}$ is a generating set for
$\cY'f =f\cY'f$. 

It follows from the description of $\cY$ via generators and relations given
after Lemma~\ref{lem:T inv} that there is a surjective $\ZZ_\ell$-algebra
homomorphism $\cY_{\psi_q}\to\cY'f$ which sends the image
$\tilde y_t\in\cY_{\psi_q}$ of $ y_t$ to
$\frac{q^{-1}\ell^a} {q-1} t' f$ and $\tilde y_{\br_i}$ to $-s_if$. By
Theorem~\ref{thm:free Cox}, $\cY_{\psi_q}$ is $\ZZ_\ell$-free of rank $|TW|$
and the same is true for $\cY'f$. The last assertion can be seen by considering
the standard basis $\{e_{U} n e_{U}\mid n \in N_G(T_0)\} $ of $\cY'$ given by
the Bruhat decomposition. Thus $\cY_{\psi_q}$ is isomorphic to $f\cY'f$ as
claimed.
\end{proof}

\subsection{Proofs of Theorem~\ref{t:main} and Corollary~\ref{cor:main}}   \label{subsec:pf thm1}
Throughout this subsection $\ell >2$ is a prime, $q$ is a prime power with
$q\equiv 1\pmod\ell$ and $a >0$ such that $\ell^a ||(q-1)$. Let $\GG=(W, L)$ be
a simply connected $\ZZ_\ell$-spets with $W$ an $\ell'$-group. Let $\cF$ be the
fusion system associated to $(W, L)$ as described in Section~\ref{subsec:B_0},
with underlying $\ell$-group $S$. Recall that with the stated assumptions we
have $S = T \cong (\ZZ/\ell^a)^n$ where $T$ is the homocyclic group
$L/{\ell^a}L$ of exponent $\ell^a$. We let $\cY$ be the Yokonuma algebra
associated to $(W, L, q)$ as in Section~\ref{sec:Yoku}.

Recall the indeterminates $\tilde u_{rj}$ with
$\tilde u_{rj}^z=\zeta_{o(r)}^{-j}u_{rj}$. By \cite[Cor.~4.8]{Ma99}, $z$ may
be chosen to divide the order of the group of roots of unity in $\QQ_\ell$,
that is, $\ell-1$. As $\ell^a||(q-1)$, by Hensel's lemma there is a unique root
of $X^z-q\in\ZZ_\ell[X]$ in $\ZZ_\ell$, say $q^{1/z}$, with
$\ell^a||(q^{1/z}-1)$. Let $\ZZ_\ell\subseteq\cO$ be a complete discrete
valuation ring containing the $|T|$th roots of unity.
Let $$\psi_\sq: \cO [\tilde\bu^{\pm1}] \to \cO, \qquad
  \tilde u_{rj}\mapsto\begin{cases} q^{\frac{1}{z}}& \text{if $j=o(r)$,}\\
   1& \text{if $1\le j<o(r)$},\end{cases}$$
be the specialisation $\psi_\sq$ from Section~\ref{subsec:spec} with $R=\cO $.
Since $|T_r|$ divides $\ell^a$, $\psi_\sq$ extends to an $\cO$-linear
homomorphism $\cO[\tilde\bu^{\pm1},\bv]\to\cO$ which we still denote $\psi_\sq$.
Let $\tilde K=\Frac(\cO[\tilde\bu^{\pm1}])$ and $K=\Frac(\cO)$.

We restate Theorem~\ref{t:main}.

\begin{thm}   \label{thm:omni}
 Let $\GG$ be as above. Suppose that $\cH(W,\bu)$ is strongly symmetric as in
 Definition~\ref{def:strongly quasi}. Then Conjectures~\ref{conj:dim B},
 \ref{conj:quot} and~\ref{conj:decomp} hold for $\GG$.
\end{thm}

\begin{proof}
Let $W_0$ be a reflection subgroup of $W$. By Lemma~\ref{lem:transfer symm},
$K\otimes_\cO\cH_{\psi_\sq}(W_0,\bu_0)$ is split semisimple and $\psi_\sq$
induces a bijection $\Irr(\tilde K \otimes_{A}\cH(W_0,\bu_0))\to
\Irr(K\otimes_\cO\cH_{\psi_\sq} (W_0,\bu_0))$. Also recall from
Section~\ref{subsec:spec}, $\Irr(\tilde K\otimes_{A} \cH(W_0, \bu_0))$ is
identified with $\Irr(W_0)$ via $\psi_1$. Henceforth we identify
$\Irr(K\otimes_\cO\cH_{\psi_\sq}(W_0,\bu_0))$ and $\Irr(W_0)$ via the bijections
induced by $\psi_\sq$ and $\psi_1$.
 
Denoting the restriction of $\psi_\sq$ to $\hA$ again by $\psi_\sq$, set
$\cY_q:= \cY_{\psi_\sq}$. Since $W$ is an $\ell'$-group, $W_\th$ is a parabolic 
subgroup of $W$ for all $\th\in\Irr_K(T)$. Then by
Theorem~\ref{thm:constr rep}(a) and the above $K\otimes_\cO \cY_q$ is split
semisimple and $\Irr(K\otimes_\cO\cY_q)$ is in bijection with pairs $(\th,\phi)$
as $\th$ runs over representatives of $W$-orbits of $\Irr(T)$ and
$\phi\in\Irr(W_\th)$. Let
$\chi'_{\th,\phi}$ be the irreducible character corresponding to the pair
$(\th,\phi)$. Then $\chi'_{\th,\phi}$ is afforded by the simple module 
$U_{\th,\phi}':=\Ind_{\hB_\th}^\hB(U)$ for $U$ a simple
$K\otimes_\cO\cH_{\psi_\sq}(W_\th,\bu_\th)$-module corresponding to $\phi$.
Here as in Section~\ref{subsec:trace form}, $\bu_\th =\bu_0$ in the notation of
Section~\ref{subsec:braid} if $W_\th =W_0$.

We consider the following $K$-linear form on $K \otimes_\cO \cY_q$: 
\begin{equation} \label{eqn:sp form}
 \tau_q:=
   \frac{1}{ |T|} \sum_{\th/W}\sum_{\phi\in\Irr(W_\th)}
   \frac{1}{\psi_\sq(f_{\th,\phi}) }\chi'_{\th,\phi}. 
\end{equation} 
By Lemma~\ref{lem:non-vanishing schur} this is well defined.
Since the coefficient of every irreducible character is non zero, $\tau_q$ is a
symmetrising form on $K\otimes_\cO \cY_\psi$ with Schur elements
$|T|\psi_\sq(f_{\th,\phi})$.
 
Let $\bW\subset B(W)$ be as in Definition~\ref{def:strongly quasi}.
By Theorem~\ref{thm:O finite generation},
$\{\tilde y_t \tilde y_\bw\mid t\in T,\ \bw\in\bW\}$ is an $\cO$-basis
of~$\cY_q$. Here, as earlier, for $x\in\cY$, we write $\tilde x:=
1_\cO\otimes x\in\cY_q$. As in Proposition~\ref{prop:is form}, we have
$$\tau_q(\tilde y_t\tilde y_\bw)=\delta_{t,1}\delta_{\bw,1}
  \qquad\text{for any $t\in T$, $\bw\in\bW$},$$
hence the above gives that the restriction of $\tau_q$ to $\cY_q$ takes values
in $\cO$. 
 
By the strongly symmetric hypothesis, and by Theorem~\ref{thm:KOW2} there
is an $\cO$-algebra isomorphism $\sigma: \cO [TW]\to \cY_q$ whose restriction to
$T$ is the identity on $T$ (where we identify $T$ with its image in $\cY_q$ via
$t\mapsto \tilde y_t$, $t\in T$). Denote also by $\sigma$ the extension
$K[TW]\to K\otimes_\cO\cY_q$. Then $\tau_\sigma:=\tau_q \circ\sigma:K[TW]\to K$
is a symmetrising form on $K[TW]$, with Schur element
$|T|\psi_\sq(f_{\th,\phi})$
at the irreducible character of $K[TW]$ corresponding under $\sigma$ to the
character $\chi'_{\th,\phi}$ of $K\otimes_\cO \cY_q$. Further, by the above
$\tau_\sigma(t) = \delta_{t,1}$ for all $t\in T$ and the
restriction of $\tau\sigma$ to $\cO [TW]$ takes values in $\cO$. Thus, by
Lemmas~\ref{lem:reversedivisibility} and~\ref{lem:otherdivisibility}, 
$$\frac{1}{|T|^2} \sum_{\th/W}\sum_{\phi\in\Irr(W_\th)}\frac{1}{\psi_\sq(f_{\th,\phi}^2)}
  =\frac{\al}{|TW|},$$
where $\al\in \cO$ is such that $\al\equiv 1\pmod\ell$. 

Recall that $S=T$ and as explained in Section~\ref{subsec:deco}, there exists
a $W$-equivariant bijection between $\Irr(T)$ and $T $. Thus the left hand side
of the above equals
$$\frac{\psi_\sq(d)} {|T|^2\psi_\sq(p_W^2 )},$$
where 
$d:= \sum_{s/\cF}\sum_{\phi\in\Irr(W(s))}{p_W^2}{f_{s,\phi }^{-2}}$.
By Lemma~\ref{lem:degree schur base}, $\psi_\sq(d) = \dim(B_0)|_{x=q}$. Since 
$\psi_\sq (p_W)\equiv |W| \pmod \ell$, we obtain the validity of
Conjectures~\ref{conj:dim B} and~\ref{conj:quot} from the displayed equation
above. 

Finally, we prove Conjecture~\ref{conj:decomp}. First of all note that since
$\sigma$ is the identity on $T$, for any pair $(\th,\phi)$ as above the
irreducible character of $K[TW]$ corresponding under $\sigma$ to the
character $\chi'_{\th,\phi}$ of $K\otimes_\cO \cY_q$ covers $\th$ and therefore
is of the form $\gamma_{\th,\tilde\phi}:=\Ind_{W_\th}^W(\tilde\phi)$ for some
$\tilde\phi\in\Irr(W_\th)$. In particular, $\sigma$ induces a permutation
$\phi\mapsto\tilde \phi$ of $\Irr(W_\th)$.
By Lemma~\ref{lem:projschur} we have that for any $\nu \in \IBr(TW)$,
$$ \sum_{\th/W}\sum_{\phi\in\Irr(W_\th)}\frac{d_{\gamma_{\th,\tilde \phi}\nu} }{\psi_\sq(f_{\th,\phi}) }$$
is divisible by $|T|$ in $\cO$. Choose a $W$-equivariant bijection between
$\Irr(T)$ and $T$ and let $\Theta:\Irr(TW)\to\Irr(B_0)$,
$\gamma\mapsto\hat\gamma$, be the bijection such that if
$\gamma = \Ind_{W_\th}^W(\tilde\phi)$, then $\hat \gamma$ is the element of
$\Irr(B_0)$ labelled by $(x,\phi)$, where the $W$-class of $x\in T$ corresponds
to the $W$-class of $\theta$ for the chosen $W$-equivariant bijection between
$\Irr(T)$ and $T$. Then $\Theta$ is $W$-equivariant. Moreover, by
Lemma~\ref{lem:degree schur base} the above displayed expression equals
$\psi_\sq(p_W) ^{-1}(\deg\Phi_{\hat \nu})|_{x =q}$. The result follows since
$\psi_\sq (p_W)$ is an invertible element of $\cO$.
\end{proof}

\begin{rem}
The above holds in a more general setting. Drop the assumption that $W$ is
spetsial; so $W$ is an $\ell$-adic reflection group of order prime to $\ell$.
Let $\psi_q:\cO [\tilde \bu^{\pm 1},\bv] \to\cO$ be any specialisation as
in Section~\ref{subsec:spec}. Suppose that $\tau:\cH(W,\bu)\to A$ is a
symmetrising form such that the following holds:
\begin{enumerate}
  \item[\rm(1)] there is a section $W\to\bW\subset B$ of the natural map
   $B\to W$ containing~$1$ whose image in $\cH(W,\bu)$ is an $A$-basis
   of $\cH(W,\bu)$ with $\tau(h_\bw) =\delta_{\bw,1}$ for all $\bw\in\bW$; and
  \item[\rm(2)] for any parabolic subgroup $W_0\le W$,
   $\tau|_{\cH(W_0,\bu_0)}:\cH(\bW_0, \bu_0)\to A$ is a symmetrising form.
\end{enumerate}
For $s\in S$, $\phi\in\Irr(W(s))$ let $f_{\tau,s,\phi}$ denote the Schur
element of $\tau|_{\cH(W_0,\bu_0)}$ with respect to $\phi$ where $W_0 = W(s)$.
Set
$$d: =\sum_{s/\cF}\sum_{\phi\in\Irr(W(s))} \frac{p_W^2}{f_{\tau,s,\phi }^2}$$
and for $\nu\in\IBr(SW)$ set 
$$\Phi_\nu(1):=
  \sum_{\gamma\in\Irr(SW)}d_{\gamma\nu}\frac{p_W}{f_{\tau, s,\phi}},$$
where $d_{\gamma\nu}$ is the decomposition number in $SW$ with respect to
$\gamma$ and $\nu$. Then with the same proof as above we get
 \begin{enumerate}
  \item[\rm(a)] $\psi_\sq(d)_\ell = |S|$;
  \item[\rm(b)] $\frac{\psi_\sq(d)}{\psi_\sq(p_W)\,|S|}\equiv1\pmod\ell$; and
  \item[\rm(c)] for each $\nu\in\IBr(SW)$, $|S|$ divides
   $\psi_\sq(\Phi_\nu(1))$.
 \end{enumerate}
\end{rem}

\begin{proof}[Proof of Corollary \ref{cor:main}]
By \cite[\S3]{MaICM} all imprimitive irreducible spetsial reflection groups are
either Coxeter groups, of type $G(e,1,n)$ or of type $G(e,e,n)$ with $n\ge3$ and
therefore strongly symmetric by Proposition~\ref{p:cases}. Therefore
Conjectures~\ref{c:main} and~\ref{c:main2} hold for these groups by
Theorem~\ref{thm:omni}. For the primitive groups, Conjecture~\ref{c:main} holds
by Proposition~\ref{prop:prim}. Conjecture~\ref{c:main2} holds when $W$ is
primitive and $2$-dimensional by Proposition~\ref{p:cases}(c) and
Theorem~\ref{thm:omni}; for $G_{14}$ Conjecture 2 holds by direct computation
using the description of the decomposition matrix provided in the proof of
Proposition~\ref{prop:Brauerdegreeconsistent}.
\end{proof}



\begin{thebibliography}{99}

\bibitem{Be01}
{\sc D. Bessis}, Zariski theorems and diagrams for braid groups. \emph{Invent.
  Math. \bf145} (2001), 487--507. 

\bibitem{BCCK}
{\sc C. Boura, E. Chavli, M. Chlouveraki, and K. Karvounis}, The BMM
  symmetrising trace conjecture for groups $G_4, G_5, G_6, G_7, G_8$.
  \emph{J. Symbolic Comput. \bf96} (2020), 62--84. 

\bibitem{BM97}
{\sc K. Bremke and G. Malle}, Reduced words and a length function for
  $G(e,1,n)$. \emph{Indag. Mathem. \bf8} (1997), 453--469.

\bibitem{BM07}
{\sc C. Broto and J. M\o ller}, Chevalley $p$-local finite groups. \emph{Algebr.
  Geom. Topol. \bf7} (2007), 1809--1919.

\bibitem{BM93}
{\sc M. Brou\'e and G. Malle}, Zyklotomische Heckealgebren. \emph{Ast\'erisque}
  No.~212 (1993), 119--189.

\bibitem{BMM99}
{\sc M. Brou\'e, G. Malle, and J. Michel}, Towards spetses I. \emph{Transform.
  Groups \bf4} (1999), 157--218.

\bibitem{BMM14}
{\sc M. Brou\'e, G. Malle, and J. Michel}, Split spetses for primitive
  reflection groups. \emph{Ast\'erisque} No.~359 (2014).

\bibitem{BMR}
{\sc M. Brou\'e, G. Malle, and R. Rouquier}, Complex reflection groups, braid
  groups, Hecke algebras. \emph{J. reine angew. Math. \bf500} (1998), 127--90.

\bibitem{Ca18}
{\sc M. Cabanes}, Local methods for blocks of finite simple groups.
  \emph{Local representation theory and simple groups}, 179--266, EMS Series
  of Lectures in Mathematics, European Mathematical Society Publishing House,
  2018.

\bibitem{Ch18}
{\sc E. Chavli}, The BMR freeness conjecture for the tetrahedral and octahedral
  families. \emph{Comm. Algebra \bf46} (2018), 386--464. 

\bibitem{CC20}
{\sc E. Chavli and M. Chouveraki}, The freeness and trace conjectures for
  parabolic Hecke algebras. Preprint. arXiv: 2007.11535.

\bibitem{ChBook}
{\sc M. Chlouveraki}, \emph{Blocks and Families for Cyclotomic Hecke Algebras}.
  Lecture Notes in Mathematics, 1981. Springer-Verlag, Berlin, 2009.

\bibitem{CdA16}
{\sc M. Chlouveraki and L. P. d'Andecy}, Markov traces on affine and cyclotomic
  Yokonuma-Hecke algebras. \emph{Int. Math. Res. Not. IMRN \bf14} (2016),
  4167--4228. 

\bibitem{CR81}
{\sc C.W. Curtis and I. Reiner}, \emph{Methods of Representation Theory.
  Vols.~1, 2.} John Wiley \& Sons, New York, 1981.

\bibitem{EGKL}
{\sc F. Eisele, M. Geline, R. Kessar, and M. Linckelmann}, Tate duality and a
  projective scalar property for symmetric algebras. \emph{Pacific J. Math.
  \bf293} (2018), 277--300.

\bibitem{GIM00}
{\sc M. Geck, L. Iancu, and G. Malle}, Weights of Markov traces and generic
  degrees. \emph{Indag. Mathem. \bf11} (2000), 379--397.

\bibitem{GP}
{\sc M. Geck and G. Pfeiffer}, \emph{Characters of Finite Coxeter Groups and
  Iwahori-Hecke Algebras.} The Clarendon Press, Oxford University Press,
  New York, 2000.

\bibitem{JK01}
{\sc J. Juyumaya and S. Kannan}, Braid relations in the Yokonuma-Hecke algebra.
  \emph{J. Algebra \bf239} (2001), 272--297.

\bibitem{KMS}
{\sc R. Kessar, G. Malle, and J. Semeraro}, Weight conjectures for
  $\ell$-compact groups and spetses. Preprint. arXiv:2003.07213, 2020.

\bibitem{Li18}
{\sc M. Linckelmann}, \emph{The Block Theory of Finite Group Algebras. Vol.~1}.
  London Math. Soc. Students Texts, {\bf 91}, Cambridge University Press,
  Cambridge, 2018.

\bibitem{Lu7}
{\sc G. Lusztig}, Character sheaves on disconnected groups. VII.
  \emph{Represent. Theory \bf9} (2005), 209--266. 

\bibitem{Ma95}
{\sc G. Malle}, Unipotente Grade imprimitiver komplexer Spiegelungsgruppen.
  \emph{J. Algebra \bf177} (1995), 768--826.

\bibitem{Ma97}
{\sc G. Malle}, Degr\'es relatifs des alg\`ebres cyclotomiques associ\'ees aux
  groupes de r\'eflexions complexes de dimension deux. \emph{Finite Reductive
  Groups (Luminy, 1994)}, 311--332, Progr. Math., 141, Birkh\"auser, Boston,
  MA, 1997.

\bibitem{MaICM}
{\sc G. Malle}, Spetses. \emph{Doc. Math. Extra Vol. ICM II} (1998), 87--96. 

\bibitem{Ma99}
{\sc G. Malle}, On the rationality and fake degrees of characters of cyclotomic
  algebras. \emph{J. Math. Sci. Univ. Tokyo \bf6} (1999), 647--677.

\bibitem{Ma00}
{\sc G. Malle}, On the generic degrees of cyclotomic algebras. \emph{Represent.
  Theory \bf4} (2000), 342--369. 

\bibitem{Ma06}
{\sc G. Malle}, Splitting fields for extended complex reflection groups and
  Hecke algebras. \emph{Transform. Groups \bf11} (2006), 195--216.

\bibitem{Ma07}
{\sc G. Malle}, Height 0 characters of finite groups of Lie type.
  \emph{Represent. Theory \bf11} (2007), 192--220.

\bibitem{MM98}
{\sc G. Malle and A. Mathas}, Symmetric cyclotomic Hecke algebras. \emph{J.
  Algebra \bf205} (1998), 275--293.

\bibitem{Ma18a}
{\sc I. Marin}, Artin groups and Yokonuma-Hecke algebras. \emph{Int. Math. Res.
  Not. IMRN \bf13} (2018), 4022--4062.

\bibitem{Ma18b}
{\sc I. Marin}, Lattice extensions of Hecke algebras. \emph{J. Algebra \bf503}
  (2018), 104--120. 

\bibitem{Ma19}
{\sc I. Marin}, Proof of the BMR conjecture for $G_{20}$ and $G_{21}$. \emph{J.
  Symbolic Comput. \bf92} (2019), 1--14. 

\bibitem{OS82}
{\sc P. Orlik and L. Solomon}, Arrangements defined by unitary reflection
  groups. \emph{Math. Ann. \bf 261} (1982), 339--357.

\bibitem{Pu94}
{\sc L. Puig}, On Joanna Scopes' criterion of equivalence for blocks of
  symmetric groups. \emph{Algebra Colloq. \bf1} (1994), 25--55.

\bibitem{Ts20}
{\sc S. Tsuchioka}, BMR freeness for icosahedral family. \emph{Exp. Math.
  \bf29} (2020), 234--245.

\bibitem{Yo67}
{\sc T. Yokonuma}, Sur la structure des anneaux de Hecke d'un groupe de
  Chevalley fini. \emph{C. R. Acad. Sci. Paris Ser. I Math. \bf264} (1967),
  344--347.

\end{thebibliography}
\end{document}